\documentclass{elsarticle}
\usepackage{amssymb}
\usepackage{amsmath}
\usepackage{hyperref}
\usepackage{bm}

\usepackage{lineno,hyperref}
\modulolinenumbers[5]

\usepackage{hyperref}
\usepackage{moreverb}

\usepackage{amsmath}
\usepackage{tikz-cd}
\usepackage[american]{babel}
\usetikzlibrary{positioning}
\tikzset{mynode/.style={draw,text width=4.4cm,align=center}
}

\newcommand{\BR}{{\mathbb R}}
\newcommand{\BC}{{\mathbb C}}
\newcommand{\BN}{{\mathbb N}}
\newcommand{\BZ}{{\mathbb Z}}
\newcommand{\Qh}{\left(-\frac{\pi}{h},\frac{\pi}{h}\right]^n}

\newcommand{\g}{\textbf{g}}

\newcommand{\RL}{{\mathbb D}_t}

\newcommand{\e}{{\bf e}}
\renewcommand{\a}{\textbf{a}}
\renewcommand{\b}{\textbf{b}}
\newcommand{\f}{{\mathfrak f}}

\newcommand{\cl}{C \kern -0.1em \ell}

\newcommand{\proof}{\bf {Proof:} \rm}

\newtheorem{theorem}{Theorem}[section]
\newtheorem{remark}{Remark}[section]
\newtheorem{lemma}{Lemma}[section]

\journal{arXiv.org}

\begin{document}

\begin{frontmatter}

\title{On fractional semidiscrete Dirac operators of L\'evy-Leblond type}
\tnotetext[mytitlenote]{\href{https://nelsonfaustino.org/}{N.~Faustino} was supported by The Center for Research and Development in Mathematics and Applications (CIDMA) through the Portuguese Foundation for Science and Technology (FCT), references UIDB/04106/2020 and UIDP/04106/2020.}

\author{N.~Faustino\corref{mycorrespondingauthor}\fnref{myfootnote}}
\address{Department of Mathematics, University of Aveiro}
\fntext[myfootnote]{ORCiD:~\href{https://orcid.org/0000-0002-9117-2021}{0000-0002-9117-2021}}

\author[mymainaddress]{Center for R\&D in Mathematics and Applications (CIDMA)}

\cortext[mycorrespondingauthor]{Corresponding author}

\address[mymainaddress]{Campus Universit\'ario de Santiago, 3810-193 Aveiro, Portugal}

\begin{abstract}
	In this paper we introduce a wide class of space-fractional and time-fractional semidiscrete Dirac operators of L\'evy-Leblond type on the semidiscrete space-time lattice $h\BZ^n\times[0,\infty)$ ($h>0$), resembling to fractional semidiscrete counterparts of the so-called parabolic Dirac operators.
	The methods adopted here are fairly operational, relying mostly on the algebraic manipulations involving Clifford algebras, discrete Fourier analysis techniques as well as standard properties of the analytic fractional semidiscrete semigroup $\left\{\exp(-te^{i\theta}(-\Delta_h)^{\alpha})\right\}_{t\geq 0}$, carrying the parameter constraints $0<\alpha\leq 1$ and $|\theta|\leq \frac{\alpha \pi}{2}$. The results obtained involve the study of Cauchy problems on $h\BZ^n\times[0,\infty)$.
\end{abstract}


\end{frontmatter}


\section{Introduction}

\subsection{The State of Art}
The study of null solutions of Dirac-like operators is known as the heart of several function theories in the context of Clifford algebras. From the fact that Dirac-like operators factorize the Laplace operator and its analogues, it is of foremost importance to study them algebraically and analytically in order to obtain refinements of well known results in harmonic analysis and to develop applications in the fields of mathematical physics and applied mathematics as well.

In 1967 L\'evy-Leblond investigates the factorization of the Schr\"odinger operator with the aim of obtaining non-relativistic analogues of the Dirac operator in the $(1+3)$--dimensional space carrying any spin (cf.~\cite{LevyLeblond67}).However, we had to wait for contemporary times to realize how the $(1+n)-$dimensional generalization of the L\'evy-Leblond type picture, coined as parabolic Dirac type operators, can be adopted to cover a wide range of applications on the crossroads of function theory and boundary value problems. Between the many notable results regarding this new class of operators, it is worth to stress the importance of the groundbreaking works of Cerejeiras, K\"ahler \& Sommen \cite{CKS05} and Cerejeiras, Sommen \& Vieira \cite{CSV07}, who have proposed the matter of factorizing the heat operator and its relatives. Since then, their work have conducted many researchers to explore this type of approach to a wide variety of model problems, including higher-dimensional analogues of the nonlinear Schr\"odinger equation \cite{B06,CFV08}.

Interestingly, the L\'evy-Leblond picture had continue to receiving an increasing interest during the last decade. Mainly, on several physical-phase-space formulations of spinning particles through supersymmetric Lie algebraic representations of PDEs (see \cite{AKTT16,A18}) as well as on operational models involving parabolic Dirac operators and their fractional analogues (see e.g. \cite{FV16,FV17,BCBM20,BRS20} and the references therein), to mention a few. Therewith, it is reasonable to say that the L\'evy-Leblond framework is no longer just an emerging topic in the mathematics and physics communities, so that nowadays one can say that its foundations and applications are well understood by many scholars.

In parallel, the study of discrete Dirac operators has experienced a rapid increase in hypercomplex analysis during the last two decades, strongly influenced by the pioneering papers \cite{faustino2006difference,FK07,FKS07,BSSV08,RSKS10}.
The literature towards this topic, whose physical roots may be found e.g. on the research papers of Kogut \& Susskind \cite{KS75} and Rabin \cite{Rabin82} (see also \cite{KK04,Sushch14} and the references therein) is very diverse. To the purpose of this paper we will depict only an abridged overview of it to motivate our approach.

For the construction of a faithful discretization for the discrete Dirac operator on the lattice $h\BZ^n$, say  $D_h$, whose precise definition will be introduced afterwards, the Clifford algebra of signature of $(n,n)$ may be embody on the ladder structure of $D_h$, with the aid of the so-called Witt basis $\left\{\e_j^+,\e_j^-~:~j=1,2,\ldots,n\right\}$, to seamlessly encode the canonical structure of the exterior algebra (cf.~\cite[Chapter 1]{gurlebeck1997quaternionic} and \cite[Chapter 2]{VR16}).

Originally highlighted in \cite{FKS07} and on the series of papers \cite{BSSV08,RSKS10}, such fact was only duly clarified in author's recent paper \cite{FaustinoKGordonDirac16}, through the canonical isomorphism $\cl_{n,n}\cong \mbox{End}(\cl_{0,n})$ between the Clifford algebra of signature $(n,n)$, $\cl_{n,n}$, and the algebra of endomorphisms acting on the Clifford algebra of signature $(0,n)$, $\cl_{0,n}$. 
This in turn allows us to establish a canonical correspondence between the discrete counterpart of the Dirac-K\"ahler operator, $d-d^*$, and the multivector approximation of Dirac operator, $D_h$, over $h\BZ^n$ (see also \cite{FaustinoMMAS17}). 

In a larger extend, 
starting from forward and backward discretizations of the Dirac operator $\displaystyle D=\sum_{j=1}^{n}\e_j\partial_{x_j}$, $D_h^+$ and $D_h^-$ respectively, formely introduced in \cite{FK07}, one can make use of the wedge ($\wedge$) and the dot $(\bullet)$ actions on $\cl_{0,n}$ to establish the canonical one-to-one correspondences
\begin{eqnarray*}
	d\overset{\text{1-1}} \longleftrightarrow D_h^- \wedge (\cdot) & \mbox{and}& d^* \overset{\text{1-1}} \longleftrightarrow D_h^+ \bullet (\cdot)
\end{eqnarray*}
in a way that the resulting discrete Dirac operator $D_{h}
$ is nothing else than a $\cl_{n,n}-$valued representation of the multivector counterpart 
of the Dirac-K\"ahler operator, $\displaystyle D_h^- \wedge (\cdot) -D_h^+ \bullet (\cdot)$ (cf.~\cite[Subsection 2.3]{FaustinoKGordonDirac16}).

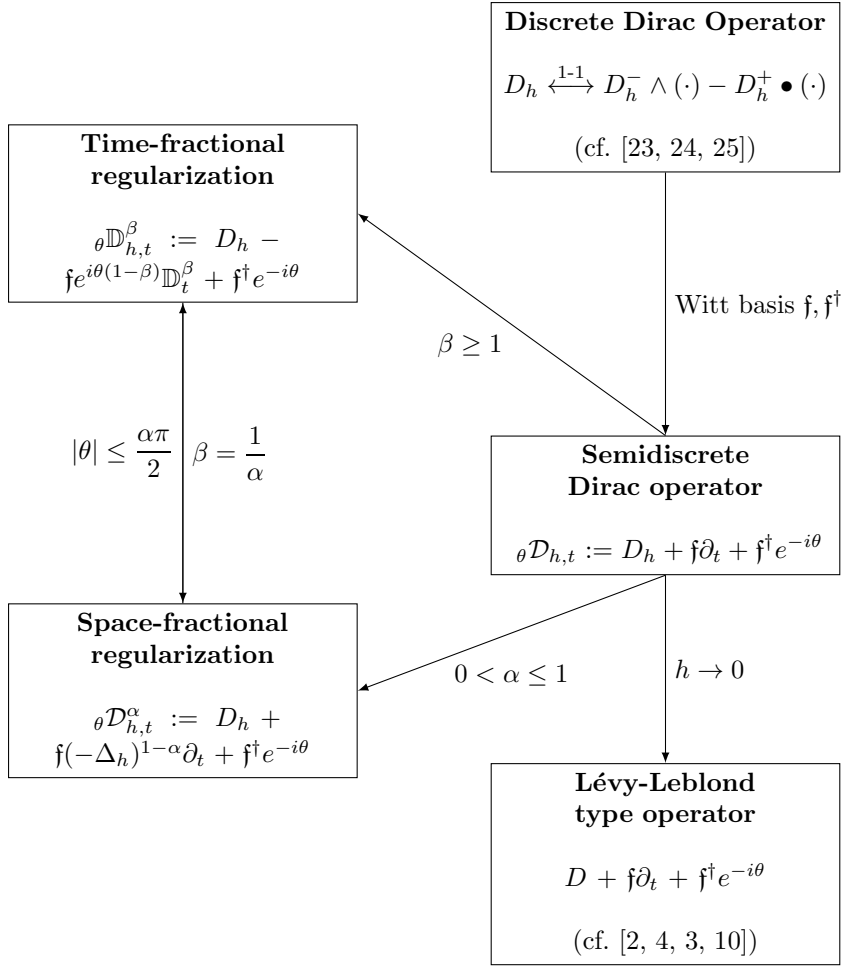
\begin{figure}\label{LevyLeblondPicture}
	\begin{tikzpicture}
		\node[mynode] (v1){\textbf{Time-fractional regularization}\\ \ \\ $	{~}_\theta\mathbb{D}_{h,t}^\beta:=D_h-\f e^{i\theta(1-\beta)} \mathbb{D}_t^\beta +\f^\dagger e^{-i\theta}$};
		\node[mynode,below right=2.5cm of v1](v2) {
			\textbf{Semidiscrete Dirac operator} \\ \ \\ ${~}_\theta\mathcal{D}_{h,t}:=\displaystyle D_h+\f \partial_t +\f^\dagger e^{-i\theta} $};
		
		\node[mynode,below = 4.0cm of v1] (v4){\textbf{Space-fractional regularization}  \\ \ \\   ${~}_\theta\mathcal{D}_{h,t}^\alpha:=D_h+ \f(-\Delta_h)^{1-\alpha}\partial_t + \f^\dagger e^{-i\theta}$};
		\node[mynode,below=2.5cm of v2](v3) {{\bf L\'evy-Leblond type operator} \\ \ \\
			$D+\f\partial_t+\f^\dagger e^{-i\theta}$ \\ \ \\ (cf.~\cite{CKS05,B06,CSV07,BCBM20})};
		\node[mynode, above=3.5cm of v2](v5) {{\bf Discrete Dirac Operator} \\ \ \\
			$D_h\overset{\text{1-1}} \longleftrightarrow D_h^-\wedge (\cdot)-D_h^+\bullet (\cdot)$  \\ \ \\
			(cf.~\cite{FaustinoKGordonDirac16,FaustinoMMAS17,F19b})};
		
		\draw[-latex] (v2.north) -- node[auto,] {$\beta \geq 1$} (v1.east);
		\draw[-latex] (v1.south) -- node[auto,] {$\beta=\dfrac{1}{\alpha}$} (v4.north);
		\draw[-latex] (v4.north) -- node[auto,] {$|\theta|\leq \dfrac{\alpha\pi}{2}$} (v1.south);
		\draw[-latex] (v2.south) -- node[below=3mm, align=center] {$0<\alpha\leq 1$} (v4.east);
		\draw[-latex] (v5.south) -- node[auto,]{Witt basis $\f,\f^\dagger$} (v2.north);
		\draw[-latex] (v2.south) -- node[auto,]{$h\rightarrow 0$} (v3.north);
	\end{tikzpicture}
	\caption{The L\'evy-Leblond picture on the fractional semidiscrete case.}
\end{figure}

\subsection{Main targets}
In this paper, we will focus our attention 
on the time-fractional and space-fractional Dirac operators,
carrying the parameters $0<\alpha\leq 1$, $\beta\geq 1$ and $|\theta|\leq\frac{\alpha \pi}{2}$. 
On the construction neatly summarized in {\bf Figure 1}, the notation $\Delta_h$ stands for the discrete Laplacian on the lattice $h\BZ^n$ considered in several author's contributions such as \cite{FaustinoKGordonDirac16,FaustinoMMAS17}.

In addition, the family of fractional discrete operators $(-\Delta_h)^{\sigma}$ ($0<\sigma\leq 1$) will be defined in terms of its Fourier symbol in the streamlines of \cite[Section 6.]{LizamaRoncal18} (see also \cite[Section 21.4.3]{F19}). For our main purposes, we will adopt the notation $\mathbb{D}_{t}^\beta$ for the so-called right-sided Riemann-Liouville fractional derivative (cf.~\cite[Chapter 2]{SKM93}).  
The operator $\mathbb{D}_{t}^\beta$, defined as follows:
\begin{eqnarray} 
	\label{RiemannLiouville}
	\mathbb{D}_{t}^\beta\Psi(x,t)=	\left\{\begin{array}{lll} 
		\displaystyle 	(-\partial_t)^k \int_{t}^{+\infty}g_{k-\beta}(s-t)\Psi(x,s)ds & \mbox{for} & k-1<\beta<k
		\\ \ \\
		(-\partial_t)^k \Psi(x,t)& \mbox{for} & \beta=k,
	\end{array}\right.
\end{eqnarray}
where $k=\lfloor \beta \rfloor+1$ ($\lfloor \beta \rfloor$ denotes the integer part of $\beta$), is an integro-differential operator for values of $\beta\neq k$, involving the higher order time-derivative $(-\partial_t)^k:=(-1)^k(\partial_t)^k$ and an integral part, corresponding to the convolution between $\Psi(x,t)$ and
the Gel'fand-Shilov function $g_{\nu}:\mathbb{R}\rightarrow [0,\infty)$, defined for $-\nu \not \in \BN_0$ by
\begin{eqnarray}
	\label{GelfandShilov}g_\nu(p)=\left\{\begin{array}{lll} 
		\displaystyle 	\frac{p^{\nu-1}}{\Gamma(\nu)} & \mbox{for} & p>0
		\\ \ \\
		0 & \mbox{for} &p\leq 0.
	\end{array}\right.
\end{eqnarray}

It should be noted that the Eulerian integral representation involving the Gamma function $\Gamma(\cdot)$ (cf.~\cite[formula 2.3.3. (1), p.322]{PBM86}):
\begin{eqnarray}
	\label{LaplaceId}\int_0^\infty e^{-p\lambda} p^{\nu-1} dp=\Gamma(\nu)\lambda^{-\nu}, & \Re(\nu)>0~~\&&~~  \Re(\lambda)>0
\end{eqnarray}
assures that for every $k-1<\beta <k$ the function $p\mapsto g_{k-\beta}(p)$, appearing on the definition of $\mathbb{D}_{t}^\beta$, defines a {\it probability density function} converging to the Dirac delta function, that is
$$
\lim_{\beta \rightarrow k^-}g_{k-\beta}(p)=\delta(p).
$$

The later set of properties allows us to say that $\mathbb{D}_{t}^\beta$ ($k-1<\beta <k$) provides us a regularization of $(-\partial_t)^k$ in the sense of the topology of the underlying space of tempered distributions.

The motivation to this paper comes from the series of contributions \cite{CiaurriGRTV18,LizamaRoncal18,GLizamaM21}, focused on the study of Cauchy problems involving fractional discrete Laplacians, and from the interrelationship between space-fractional and time-fractional operators highlighted in \cite{CLRV15} towards superdiffusion equations. Such link has played a crucial role in the theory of PDEs (cf.~\cite{KL02,KL12}) and on stochastic calculus as well (cf.~\cite{BMN09}). For a survey on both topics, we also refer to the monograph \cite{MS12}.

Thus, we are not only interested on fractional difference analogues for the parabolic operator of heat type ($\alpha=\beta=1$, $\theta=0$ \& $h\rightarrow 0$) and of Schr\"odinger type ($\alpha=\beta=1$, $\theta=\pm\frac{\pi}{2}$ \& $h\rightarrow 0$), considered in the series of papers \cite{CKS05,B06,CSV07,AKTT16,A18,BCBM20}, neither fractional counterparts of semidiscrete models involving the {\it semidiscrete heat operator} $\partial_t-\Delta_h=(D_h+\f\partial_t+\f^\dagger)^2$ (i.e when $\alpha=\beta=1$ \& $\theta=0$), already considered in~\cite{BBRS14},
but also on the interface between space-fractional and time-fractional semidiscrete operators in a way that the {\it fractional discrete Laplacian} $-(-\Delta_h)^\alpha$, carrying the parameter $0<\alpha\leq 1$, turns out be treated as a temporal regularization of order $\beta:=\frac{1}{\alpha}\geq 1$, encoded on the time-fractional derivative $\mathbb{D}_t^\beta$. 

Moreover, the parameter condition $|\theta|\leq \frac{\alpha\pi}{2}$ encoded on the [fractional semidiscrete] analytic semigroup $\displaystyle \left\{\exp(-te^{i\theta}(-\Delta_h)^\alpha)\right\}_{t\geq 0}$ is ubiquitous on space-fractional diffusion models (cf.~\cite{Mainardi01}). Its incorporation on our model problem is seamlessly justified by the constraint $|\theta|=|\arg(te^{i\theta})|<\pi$, carrying the integral representation obtained in \cite[p. 458, eq.~(21.35)]{F19} for {\it modified Bessel functions of the first kind} (cf.~\cite[Subsection 2.1 \& Subsection 3.1]{LizamaRoncal18}). This will be the main novelty of this paper in comparison with the semidiscrete heat semigroup representations considered in the series of author's contributions \cite{F19,F19b,F20}.

Up to author's knowledge, the overlap between space-fractional Dirac operators -- such as the ones introduced by Bernstein in \cite{B16} -- and time-fractional Dirac operators of Riemann-Liouville type -- such as the ones considered e.g. by Ferreira \& Vieira \cite{FV16} --  was not yet addressed so that the idea of connecting fractional order in time and space on this paper shall be seen eventually as another step further to pursue the goal of studying the mild solutions for time-fractional Navier-Stokes equations in the superdiffusive case, from a hypercomplex analysis perspective (see e.g. \cite{CarvalhoPlanas15,ZP17} and the references therein for an overview on the subdiffusive case).


\subsection{Layout of the paper}
This paper is organized as follows:
\begin{itemize}
	\item In Section \ref{DefsSection} we briefly provide some background on Clifford algebras and on {\it discrete Fourier analysis} required throughout the paper. From the point of view of the theory of pseudo-differential operators, that will allows us to shift all the well-known constructions in the space of square-integrable functions to the space of Clifford-valued distributions over the lattice $h\BZ^n$ (see, for instance, \cite[Subsection 21.2]{F19} \& \cite[Subsection 2.2.]{F19}).
	\item In Section \ref{SemidiscreteSection} we introduce a time-fractional and a space-fractional variant of the {\it semidiscrete Dirac operator} $D_h+\f\partial_t+\f^\dagger e^{-i\theta}$. 
	Guided by the approaches of Cerejeiras, K\"ahler \& Sommen \cite{CKS05} and Cerejeiras, Sommen \& Vieira \cite{CSV07}, we introduce the time-fractional regularization of $D_h+\f\partial_t+\f^\dagger e^{-i\theta}$ by replacing the time-derivative $\partial_t$ by a time-fractional counterpart $-e^{i\theta(1-\beta)}\mathbb{D}_t^\beta$ ($\beta\geq 1$), mixing the Riemann-Liouville derivative and a unitary term lying on the unit circle $\mathbb{S}^1$. For the space-fractional regularization, we consider the space-fractional counterpart $(-\Delta_h)^{\alpha-1}\partial_t$ of $\partial_t$, involving the fractional discrete operator $(-\Delta_h)^{\alpha-1}$ ($0<\alpha\leq 1$) instead of $\partial_t$ (case of $\alpha=1$). With the proof of {\bf Theorem \ref{CoupledSystemParabolicDiracRLt}} \& {\bf Theorem \ref{CoupledSystemParabolicDiract}}, we will show that the formulations highlighted on {\bf Figure 1} retain all of the salient features of the null solutions of the {\it L\'evy-Leblond type operator} $D+\f \partial_t+\f^\dagger e^{-i\theta}$ (limit case $h\rightarrow 0$) considered by several authors (see also the papers of Bernstein \cite{B06}, and Bao, Constales, De Bie \& Mertens \cite{BCBM20}).
	\item In Section \ref{MainResultsSection} we will show that the null solutions of both fractional semidiscrete operators are indeed interrelated.
	Starting from the pseudo-differential representation of the fractional semidiscrete analytic semigroup \newline $\displaystyle \left\{\exp(-te^{i\theta}(-\Delta_h)^\alpha)\right\}_{t\geq 0}$ in terms of its Fourier symbol we will provide first, with the proof of {\bf Theorem \ref{SemigroupCauchySolver}}, a bridge result that establishes the correspondence between Cauchy problems of space-fractional and time-fractional type, encoded by the set of fractional operators, $e^{-i\theta}\partial_t+(-\Delta_h)^\alpha$ and $-e^{-i\theta \beta}\mathbb{D}_t^\beta-\Delta_h$ respectively. This theorem is essentially a wise reformulation of \cite[Theorem 3]{CLRV15} in terms of the right-sided Riemann-Liouville fractional derivative.
	Afterwards, we will prove {\bf Theorem \ref{LevyLeblondSolutions}} -- the main contribution of this paper -- by combining {\bf Theorem \ref{CoupledSystemParabolicDiracRLt}}, {\bf Theorem \ref{CoupledSystemParabolicDiract}} \& {\bf Theorem \ref{SemigroupCauchySolver}}. 
	\item In Section \ref{SpaceFractionalRemarks} we will briefly discuss an alternative formulation for the {\it semidiscrete Dirac operator of space-fractional type} that factorizes the space-fractional operator $e^{-i\theta}\partial_t+(-\Delta_h)^\alpha$. In the end we will also comment on the key ingredients considered to obtain our main results and depict further directions of research on the crosssroads of function spaces and Helmholtz-Leray type decompositions.
\end{itemize}

\section{Definitions}\label{DefsSection}

\subsection{Clifford algebra setup}

Let $\e_0,\e_1,\e_2,\ldots,\e_n,\e_{n+1},\e_{n+2}\,\ldots,\e_{2n},\e_{2n+1}$ be the generators of the Clifford algebra of signature $(n+1,n+1)$, $\cl_{n+1,n+1}$, satisfying 
\begin{eqnarray}
	\label{CliffordBasis}
	\begin{array}{lll}
		\e_j \e_k+ \e_k \e_j=-2\delta_{jk}, & 0\leq j,k\leq n \\
		\e_{j} \e_{n+k}+ \e_{n+k} \e_{j}=0, & 0\leq j\leq n ~~\&~~  1\leq k\leq n+1\\
		\e_{n+j} \e_{n+k}+ \e_{n+k} \e_{n+j}=2\delta_{jk}, & 1\leq j,k\leq
		n+1.
	\end{array}
\end{eqnarray}

As it is well-known from the literature (cf.
\cite[Chapter 3]{VR16}), $\cl_{n+1,n+1}$ is a
universal algebra of dimension $2^{2n+2}$ linear isomorphic to the
exterior algebra $\Lambda^*\left(\BR^{n+1,n+1}\right)$,   containing the field of real numbers $\BR$, the $n-$dimensional Euclidean space $\BR^n$ and the Minkowski space $\BR^{n+1,n+1}$ of signature $(n+1,n+1)$ as proper subspaces. 

In particular, the ladder structure of $\cl_{n+1,n+1}$ allows us to represent the space-time tuple $(x_1,x_2,\ldots,x_n,t)$ of $\BR^{n+1}$ through the paravector representation $\displaystyle t+x=t+\sum_{j=1}^{n}x_j\e_j$ of  $\BR\oplus \BR^n$. 

Here and elsewhere, the $1-$vector
representations $\displaystyle x=\sum_{j=1}^{n}x_j\e_j$ and $\displaystyle x\pm h
\e_j$ of $\BR^{n}$ will be adopted to describe
the lattice point $(x_1,x_2,\ldots,x_n)$ of $h\BZ^{n}$ and the
forward/backward shifts $(x_1,x_2,\ldots,x_j\pm h,\ldots, x_n)$ over
$h\BZ^n$, respectively. 
For a sake of readibility, one will use throughout the manuscript the tuple notations $(x+h\e_j,t)$ and $(x-h\e_j,t)$
to define shifts over the semidiscrete space-time lattice
$$h\BZ^n \times \left[0,\infty\right):=\left\{(x,t)\in \BR^{n+1}~:~\frac{x}{h}\in \BZ^n~~\wedge~~t\geq 0 \right\}.$$ 

We notice also that the Clifford algebras $\cl_{0,n}$, $\cl_{1,1}$ and $\cl_{1,n+1}$, considered e.g. on the series of papers \cite{CKS05,B06,FV16,FV17,BCBM20}, are subalgebras of $\cl_{n+1,n+1}$. For our purposes, we assume that $\cl_{1,1}$ is generated by the nilpotents
\begin{eqnarray}
	\label{WittBasis}	\f=\frac{1}{2}(\e_{2n+1}+\e_0) & \mbox{and} & \f^\dagger=\frac{1}{2}(\e_{2n+1}-\e_0).
\end{eqnarray} 

Thereby, the graded anti-commuting relations (\ref{CliffordBasis}) are equivalent to
\begin{eqnarray}
	\label{CliffordWittBasis}
	\begin{array}{lll}
		\e_j  \f+\f \e_j=0, & \e_{n+j} \f+\f \e_{n+j}=0, &(1\leq j\leq n)\\ \ \\
		\e_j  \f^\dagger+\f^\dagger \e_j=0, & \e_{n+j} \f^\dagger +\f^\dagger \e_{n+j}=0, &(1\leq j\leq n)\\ \ \\
		(\f)^2=(\f^\dagger)^2=0,& \f \f^\dagger +\f^\dagger \f=1. &
	\end{array}
\end{eqnarray}

Furthermore, based on the decomposition $\cl_{n+1,n+1}=\cl_{1,1}\otimes \cl_{n,n}$, we will represent the Clifford-vector-valued functions $(x,t)\mapsto \Psi(x,t)$ with membership in the {\it complexified Clifford algebra} $\BC\otimes\cl_{n+1,n+1}$, through the ansatz
\begin{eqnarray*}
	\Psi(x,t)=\Psi^{[0]}(x,t)+\f\Psi^{[1]}(x,t)+\f^\dagger \Psi^{[2]}(x,t)+\f\f^\dagger \Psi^{[3]}(x,t),
\end{eqnarray*}
whereby $(x,t)\mapsto\Psi^{[m]}(x,t)$ ($m=0,1,2,3$) are Clifford-vector-valued functions with membership in $\BC \otimes \cl_{n,n}$.

To introduce in the sequel function spaces and operators underlying to the $\BC \otimes \cl_{n+1,n+1}-$valued functions $(x,t)\mapsto \Psi(x,t)$ and to its $\BC\otimes\cl_{n,n}-$valued components $(x,t)\mapsto\Psi^{[m]}(x,t)$ ($m=0,1,2,3$), as well as $\BC\otimes\cl_{n,n}-$valued representations of the discrete Dirac operator $D_h$ in terms of its Fourier multipliers (cf.~\cite[Subsection 21.2.2]{F19} and \cite[Subsection 2.3]{F20}), one 
has to consider the $\dag-${\it
	conjugation} operation ${\bf a} \mapsto{\bf a}^\dag$ on the \textit{complexified Clifford algebra} $\BC\otimes\cl_{n+1,n+1}$, defined recursively as follows: 
\begin{eqnarray}
	\label{dagconjugation}
	\begin{array}{lll}
		(\a \b)^\dag=\b^\dag\a^\dag \\ (a_J \e_J)^\dag =\overline{a_J}~\e_{j_r}^\dag
		\ldots \e_{j_2}^\dag\e_{j_1}^\dag~~~(0\leq j_1<j_2<\ldots<j_r\leq 2n+1) \\
		\e_j^\dag=-\e_j~~~\mbox{and}~~~\e_{n+1+j}^\dag=\e_{n+1+j}~~~(0\leq j\leq
		n).
	\end{array}.
\end{eqnarray}

From (\ref{dagconjugation}), the $\dagger-$conjugation identities
\begin{center}
	$(\f^\dag)^\dag=\f$ and $(\f \f^\dagger)^\dagger=\f \f^\dagger$,
\end{center}
involving the nilpotents $\f$ and $\f^\dagger$ defined viz eq. (\ref{WittBasis}), are then immediate. Also, from (\ref{dagconjugation}) one readily has that $$\a\mapsto \| \a\|:=\sqrt{\a^\dagger \a}$$
defines a $\|\cdot \|-$norm endowed by the \textit{complexified Clifford algebra} structure of $\BC\otimes \cl_{n+1,n+1}$, since $\a^\dag\a=\a\a^\dag$ is a non-negative real number. In case where $\a$ belongs to $\BC\otimes \BR^{n+1,n+1}$, the quantity $\|\a\|$ equals to the standard norm of $\a$ on $\BC^{2n+2}$.

\subsection{The discrete Fourier analysis background}

Let us define by $\mathcal{S}(h\BZ^n;\BC\otimes \cl_{n,n}):=\mathcal{S}(h\BZ^n)\otimes \left(\BC\otimes \cl_{n,n}\right)$ the Schwartz class of $\BC\otimes\cl_{n,n}-$valued functions on the lattice $h\BZ^n$, consisting on \textit{rapidly decaying functions} $x\mapsto \varPhi(x,t)$ ($t \in [0,\infty)$) defined for any $M\in [0,\infty)$ by the semi-norm condition
$$ \displaystyle \sup_{x \in h\BZ^n} (1+\| x\|^2)^M~\| \varPhi(x,t)\|<\infty,$$
and by $\ell_2(h\BZ^n;\BC\otimes \cl_{n,n}):=\ell_2(h\BZ^n)\otimes \left(\BC\otimes \cl_{n,n}\right)$ the \textit{right Hilbert module} endowed by the Clifford-valued sesquilinear form 
\begin{eqnarray}
	\label{lpInner} \langle \varPhi(\cdot,t),\varPsi(\cdot,t) \rangle_{h}=\sum_{x\in h\BZ^n}h^n~
	\varPhi(x,t)^\dag\varPsi(x,t).
\end{eqnarray}

By exploiting \cite[Exercise 3.1.7]{RuzhanskyT10} to the Clifford-valued setting, it is easy to check that the seminorm condition
$$ \displaystyle \sup_{x \in h\BZ^n} (1+\| x\|^2)^{-M}~\| \varPhi(x,t)\|<\infty$$
induces the set of all \textit{continuous linear functionals} with membership in the Schwarz class $\mathcal{S}(h\BZ^n;\BC\otimes \cl_{n,n})$, induced by the mapping $$\varPhi(\cdot,t) \mapsto \langle \varPhi(\cdot,t),\varPsi(\cdot,t)\rangle_{h},$$ whereby the family of functions $\varPsi(\cdot,t):h\BZ^n \rightarrow \BC\otimes \cl_{n,n}$ ($t\in [0,\infty)$) belong to the space of $\BC\otimes \cl_{n,n}-$valued tempered distributions on the lattice $h\BZ^n$, denoted as $$\mathcal{S}'(h\BZ^n;\BC\otimes \cl_{n,n}):=\mathcal{S}'(h\BZ^n)\otimes\left(\BC\otimes \cl_{n,n}\right).$$  

In particular, the mapping property $\varPhi(\cdot,t) \mapsto \langle \varPhi(\cdot,t),\varPsi(\cdot,t)\rangle_{h}$ together with density arguments allows us to define, for every $x\mapsto \varPhi(x,t)$ with membership in
$
\ell_2(h\BZ^n;\BC\otimes \cl_{n,n})
$, a distribution $\varPhi(\cdot,t) \mapsto \langle \varPhi(\cdot,t),\varPsi(\cdot,t)\rangle_{h}$ lying to $\mathcal{S}'(h\BZ^n;\BC\otimes \cl_{n,n})$.

Next, we denote by $\Qh$ the
$n-$dimensional \textit{Brioullin zone} representation for the $n-$torus $\BR^n/\frac{2\pi}{h}\BZ^n$ (cf.~\cite[p.~
324]{Rabin82}), by $$L_2\left(\Qh;\BC\otimes\cl_{n,n}\right):=L_2\left(\Qh\right)\otimes\left(\BC\otimes \cl_{n,n}\right)$$ the $\BC \otimes \cl_{n,n}-$\textit{Hilbert module} endowed by the sesquilinear form
\begin{eqnarray}
	\label{BilinearFormQh}	\displaystyle \langle {\bf f}(\cdot,t),\g(\cdot,t
	)\rangle_{\Qh}= \int_{\Qh} {\bf f}(\xi,t)^\dag \g(\xi,t) d\xi,
\end{eqnarray}
and by $C^\infty\left(\Qh;\BC \otimes \cl_{n,n}\right)$ the space of $\BC \otimes \cl_{n,n}-$valued test functions.
The \textit{discrete Fourier transform}, defined by
\begin{eqnarray}
	\label{discreteFh}
	(\mathcal{F}_{h} \varPhi)(\xi,t)=\left\{\begin{array}{lll}
		\displaystyle \frac{h^n}{\left(2\pi\right)^{\frac{n}{2}}}\displaystyle 
		\sum_{x\in h\BZ^n}\varPhi(x,t)e^{i x \cdot \xi}, & \xi\in \Qh
		\\ \ \\
		0  & \xi\in \BR^n \setminus \Qh
	\end{array}\right.
\end{eqnarray}
yields the isometric isomorphism $$\mathcal{F}_{h}:\ell_2(h\BZ^n;\BC \otimes \cl_{n,n})\rightarrow L_2\left(\Qh;\BC \otimes \cl_{n,n}\right),$$ whose inverse $(\mathcal{F}_{h}^{-1} \g)(x,t)=\widehat{\g}_{h}(x,t)
$ is provided by the Fourier coefficients 
\begin{eqnarray}
	\label{FourierInversion}
	\widehat{\g}_{h}(x,t)=\frac{1}{(2\pi)^{\frac{n}{2}}}\int_{\Qh} (\mathcal{F}_{h} \g)(\xi,t) e^{-i x \cdot \xi} d\xi.
\end{eqnarray} 

Moreover, by noticicing that the sesquilinear form (\ref{BilinearFormQh}) allows us to define a mapping that identifies $C^\infty\left(\Qh;\BC\otimes \cl_{n,n}\right)$ with the dual space $C^\infty\left(\Qh;\BC \otimes \cl_{n,n}\right)'$, the so-called space of $\BC\otimes \cl_{n,n}-$valued distributions over $\Qh$ (cf.~\cite[Exercise 3.1.15.]{RuzhanskyT10} \&  \cite[Definition 3.1.25]{RuzhanskyT10}),
we immediately get that the \textit{Parseval type relation}
\begin{eqnarray}
	\label{Parseval}	\langle \mathcal{F}_{h} \varPhi(\cdot,t),\g(\cdot,t
	)\rangle_{\Qh}=\left\langle \varPhi(\cdot,t),\widehat{\g}_{h}(\cdot,t)\right\rangle_{h},
\end{eqnarray}
involving the sesquilinear forms (\ref{lpInner}) and (\ref{BilinearFormQh}) (cf.~\cite[Definition 3.1.27]{RuzhanskyT10}), allows us to extend propertly $\mathcal{F}_{h}$ (see eq.~ (\ref{discreteFh})) to the setting of distributions, through the mapping property (cf.~\cite[Definition 3.1.27 \& 3.1.28]{RuzhanskyT10})
$$\mathcal{F}_{h}:\mathcal{S}(h\BZ^n;\BC\otimes \cl_{n,n})\rightarrow C^\infty\left(\Qh;\BC \otimes \cl_{n,n}\right),$$
in a way that the Fourier coefficients $x\mapsto \widehat{\g}_{h}(x,t)$, defined {\it viz} eq. (\ref{FourierInversion}), belong to $\mathcal{S}'(h\BZ^n;\BC\otimes \cl_{n,n})$.

In the same order of ideas of \cite[Subsection 21.1.3]{F19} \& \cite[Subsection 2.2]{F20} (see also \cite[Section 6]{CiaurriGRTV17}), one can also define the \textit{discrete convolution operation} $\star_{h}$ between the \textit{discrete distribution} $\varPsi(\cdot,t)\in \mathcal{S}'(h\BZ^n;\BC \otimes \cl_{n,n})$, and the \textit{discrete function} $\Phi(x)\in \mathcal{S}(h\BZ^n;\BC \otimes \cl_{n,n})$ as follows:
\begin{eqnarray}
	\label{discreteConvolution} \left(\varPsi(\cdot,t)\star_{h} \Phi\right)(x)=\sum_{y\in h\BZ^n} h^n~ \Phi(y)\varPsi(x-y,t).
\end{eqnarray}

Indeed, from the duality condition
\begin{eqnarray*}
	\left\langle~ \varPsi(\cdot,t) \star_{h} \Phi,\g(\cdot,t)~\right\rangle_{h}=\langle~ \varPsi(\cdot,t),\widetilde{\Phi} \star_{h} \g(\cdot,t)~\rangle_{h}, & \mbox{with}& \widetilde{\Phi}(x)=[\Phi(-x)]^\dagger
\end{eqnarray*}
that yields straightforwardly from the following sequence of identities:
\begin{eqnarray*}
	\langle~\mathcal{F}_{h}\left[\varPsi(\cdot,t)\star_{h} \Phi\right],\g(\cdot,t) ~\rangle_{\Qh}&=&
	\langle~\varPsi(\cdot,t)\star_{h} \Phi, \mathcal{F}_{h}^{-1}[\g(\cdot,t)]~ \rangle_{h}
	\\
	&=&\langle~\varPsi(\cdot,t), \widetilde{\Phi}\star_{h}\mathcal{F}_{h}^{-1}[\g(\cdot,t)]~ \rangle_{h} \\
	&=&\left\langle~\varPsi(\cdot,t),\mathcal{F}_{h}^{-1}\left( \mathcal{F}_{h}\widetilde{\Phi}~\g(\cdot,t)\right)~ \right\rangle_{h} \\
	&=&\left\langle~\mathcal{F}_{h}\varPsi(\cdot,t), \mathcal{F}_{h}\widetilde{\Phi}~\g(\cdot,t)~ \right\rangle_{h} \\
	&=&\langle~\left(\mathcal{F}_{h}\varPsi(\cdot,t)\right)\left(\mathcal{F}_{h}\Phi\right),\g(\cdot,t) ~\rangle_{\Qh},
\end{eqnarray*}
one can say the \textit{discrete convolution operation} (\ref{discreteConvolution}) is well-defined.

\section{Fractional semidiscrete Dirac operators of L\'evy-Leblond type}\label{SemidiscreteSection}
\subsection{The time-fractional case}\label{TimeFractionalSection}

Let us now recall the basic setup and results from the series of papers \cite{FaustinoKGordonDirac16,FaustinoMMAS17,F19,F19b} to discuss further aspects of our
construction. On the sequel, we will use the nilpotents $\f$ and $\f^\dagger$ of $\cl_{1,1}$
to introduce firstly on $\BC\otimes\cl_{n+1,n+1}$, with $\cl_{n+1,n+1}=\cl_{1,1}\otimes \cl_{n,n}$, the semidiscrete Dirac type operator 
\begin{eqnarray}
	\label{SemidiscreteDh}
	{~}_\theta\mathcal{D}_{h,t}:=D_h+\f \partial_t +\f^\dagger e^{-i\theta},
\end{eqnarray}
carrying the discrete Dirac operator
\begin{eqnarray}
	\label{DiracEqh}
	\begin{array}{lll}
		D_h \Psi(x,t)&=&\displaystyle \sum_{j=1}^n\e_j\frac{\Psi(x+h \e_j,t)-\Psi(x-h\e_j,t)}{2h}+ \\
		&+&	\displaystyle \sum_{j=1}^n\e_{n+j}\frac{2\Psi(x,t)-\Psi(x+h \e_j,t)-\Psi(x-h \e_j,t)}{2h},
	\end{array}
\end{eqnarray}
the time derivative $\partial_t$ and the unitary group parameter $\theta \mapsto e^{-i\theta}$.

From the factorization property $\left(D_h\right)^2=-\Delta_h$ involving the \textit{discrete Laplacian} (cf.~\cite[Proposition 2.1]{FaustinoKGordonDirac16})
\begin{eqnarray}
	\label{discreteLaplacian}
	\displaystyle \Delta_h \varPsi(x,t)=\sum_{j=1}^n
	\frac{\varPsi(x+h\e_j,t)+\varPsi(x-h\e_j,t)-2\varPsi(x,t)}{h^2},
\end{eqnarray}
and from the graded anti-commuting relations (\ref{CliffordWittBasis}) we can easily prove that 
\begin{eqnarray*}
	\left({~}_\theta\mathcal{D}_{h,t}\right)^2= e^{-i\theta}\partial_t-\Delta_h.
\end{eqnarray*}

Here, we notice that $\left({~}_\theta\mathcal{D}_{h,t}\right)^2$ is, up to the unitary term $e^{-i\theta}$, the semidiscrete heat operator considered in \cite{BBRS14}. 

Accordingly to \cite{FV16,FV17,BRS20}, a time-fractional counterpart of (\ref{SemidiscreteDh}) can be straightforwardly obtained by replacing the time-derivative by a fractional analogue. In particular, based on the Riemann-Liouville derivative (\ref{RiemannLiouville}) we introduce the time-fractional regularization of ${~}_\theta\mathcal{D}_{h,t}$, defined for $\beta\geq 1$ by the {\it time-fractional semidiscrete} operator
\begin{eqnarray}
	\label{TimeFractionalDh}
	{~}_\theta\mathbb{D}_{h,t}^\beta:=D_h-\f e^{i\theta(1-\beta)}\mathbb{D}_t^\beta +\f^\dagger e^{-i\theta}.
\end{eqnarray}

In addition, from the set of graded anti-commuting relations (\ref{CliffordWittBasis}) one can also obtain the factorization property 
\begin{eqnarray*}
	\left({~}_\theta\mathbb{D}_{h,t}^\beta\right)^2=-e^{-i\theta\beta}\mathbb{D}_t^\beta-\Delta_h
\end{eqnarray*}
for (\ref{TimeFractionalDh}) in terms of the {\it time-fractional semidiscrete} operator $-e^{-i\theta\beta}\mathbb{D}_t^\beta-\Delta_h$.

\subsection{Time-fractional case vs. space-fractional case}\label{SpaceFractionalSection}

By taking into account the \textit{discrete Fourier transform} introduced viz eq.~(\ref{discreteFh}), one can obtain an equivalent formulation of  ${~}_\theta\mathbb{D}_{h,t}^\beta$ resp.  $-e^{-i\theta\beta}\partial_t-\Delta_h$ in terms of the Fourier multipliers (cf.~\cite[Section 3.3]{F19b})
\begin{eqnarray}
	\label{FourierMultipliers}
	\begin{array}{lll}
		d_h(\xi)^2&=&\displaystyle \sum_{j=1}^{n}\frac{4}{h^2}\sin^2\left(\frac{h\xi_j}{2}\right) \\ \ \\ 
		\textbf{z}_{h}(\xi)
		&=&\displaystyle \sum_{j=1}^n -i\e_j\dfrac{\sin(h\xi_j)}{h}+\sum_{j=1}^n \e_{n+j}\dfrac{1-\cos(h\xi_j)}{h}
	\end{array}
\end{eqnarray} 
of $\mathcal{F}_h\circ (-\Delta_h)\circ \mathcal{F}_h^{-1}$ and $\mathcal{F}_h\circ D_h\circ \mathcal{F}_h^{-1}$. Here, we emphasize that the Fourier multiplier formulation lead to the following mapping properties (cf.~\cite[Subsection 21.2.2]{F19})
\begin{eqnarray*}
	D_h:\mathcal{S}(h\BZ^n;\BC \otimes \cl_{n,n})\rightarrow \mathcal{S}(h\BZ^n;\BC \otimes \cl_{n,n}), \\
	-\Delta_h:\mathcal{S}(h\BZ^n;\BC \otimes \cl_{n,n})\rightarrow \mathcal{S}(h\BZ^n;\BC \otimes \cl_{n,n})
\end{eqnarray*}
in a way that
\begin{eqnarray*}
	\mathcal{F}_h\left[\left(-e^{-i\theta\beta}\mathbb{D}_t^\beta-\Delta_h\right)\Psi(x,t)\right]&=&\left(-e^{-i\theta\beta}\mathbb{D}_t^\beta+d_h(\xi)^2\right)\mathcal{F}_h\Psi(\xi,t),
	\\ \ \\
	\mathcal{F}_h\left[{~}_\theta\mathbb{D}_{h,t}^\beta \Psi(x,t)\right]&=&	\left(\textbf{z}_{h}(\xi)-\f e^{i\theta(1-\beta)}\mathbb{D}_t^\beta+\f^\dagger e^{-i\theta}\right)\mathcal{F}_h\Psi(\xi,t).
\end{eqnarray*}

Thereby, by mimicking \cite[Theorem 4.1 \& Corollary 4.2]{CSV07} we obtain the following theorem.
\begin{theorem}\label{CoupledSystemParabolicDiracRLt}
	For every $t\in [0,\infty)$, let us assume that the components $\Psi^{[m]}(x,t)$ of the
	$\BC\otimes\cl_{n+1,n+1}-$valued function \begin{eqnarray*}
		\Psi(x,t)=\Psi^{[0]}(x,t)+\f\Psi^{[1]}(x,t)+\f^\dagger\Psi^{[2]}(x,t)+\f\f^\dagger\Psi^{[3]}(x,t)
	\end{eqnarray*}
	satisfy the set of conditions ($m=0,1,2,3$)
	\begin{eqnarray*}
		\Psi^{[m]}(\cdot,t)\in \mathcal{S}(h\BZ^n;\BC \otimes \cl_{n,n}) &\&& \mathbb{D}_t^\beta \Psi^{[m]}(\cdot,t)\in \mathcal{S}(h\BZ^n;\BC \otimes \cl_{n,n}). 
	\end{eqnarray*}
	
	Then, for every $\beta\geq 1$ the function $\Psi(x,t)$ is a null solution of ${~}_\theta\mathbb{D}_{h,t}^{\beta}$ (see eq.~(\ref{TimeFractionalDh})) if, and only if
	\begin{eqnarray}
		\label{CoupledSystemsHeatAlphaRL} \left\{\begin{array}{lll} 
			\RL^\beta\Psi^{[m]}(x,t)=-e^{i\theta\beta}\Delta_h\Psi^{[m]}(x,t) & \mbox{for} & m=0,2
			\\ \ \\
			\Psi^{[1]}(x,t)=-e^{i\theta}D_h\Psi^{[0]}(x,t) & &\\ \ \\
			\Psi^{[3]}(x,t)=e^{i\theta}D_h\Psi^{[2]}(x,t)-\Psi^{[0]}(x,t). &  & 
		\end{array}\right.
	\end{eqnarray}
\end{theorem}

\begin{remark}
	Under the conditions of {\bf Theorem \ref{CoupledSystemParabolicDiracRLt}}, one can say that the $\BC\otimes\cl_{n,n}-$valued components $\Psi^{[0]}(x,t)$ and $\Psi^{[2]}(x,t)$ of $\Psi(x,t)$ are null solutions of the {\it time-fractional semidiscrete operator} $-e^{-i\theta\beta}\RL^\beta-\Delta_h$.
\end{remark}

In the space-fractional case, one has to consider the \textit{discrete Fourier transform} $\mathcal{F}_h$ and its inverse $\mathcal{F}_h^{-1}$, defined viz eq.~(\ref{FourierInversion}), to introduce a space-fractional regularization of the semidiscrete Dirac operator ${~}_\theta\mathcal{D}_{h,t}$. Namely, with the aid of the fractional discrete operator $(-\Delta_h)^{\sigma}=\mathcal{F}_h^{-1}\circ \left(d_h(\xi)^2\right)^\sigma \circ \mathcal{F}_h$:
\begin{eqnarray}
	\label{fractionalDeltah}(-\Delta_h)^{\sigma}\Psi(x,t)&=&\frac{1}{(2\pi)^{\frac{n}{2}}}\int_{\Qh} \left(d_h(\xi)^2\right)^\sigma \mathcal{F}_h\Psi(\xi,t)e^{-ix\cdot \xi}d\xi,
\end{eqnarray} 
one can define
\begin{eqnarray}
	\label{SpaceFractionalDh}{~}_\theta\mathcal{D}_{h,t}^\alpha :=D_h+ \f(-\Delta_h)^{1-\alpha}\partial_t +e^{-i\theta} \f^\dagger
\end{eqnarray}
as the space-fractional regularization underlying to (\ref{SemidiscreteDh}).

We note already that the factorization property
$$({~}_\theta\mathcal{D}_{h,t}^\alpha)^2=e^{-i\theta}(-\Delta_h)^{1-\alpha}\partial_t-\Delta_h$$
yields from the set of identities
involving the Fourier multipliers of
\begin{eqnarray*}
	\mathcal{F}_h\circ {~}_\theta\mathcal{D}_{h,t}^\alpha \circ \mathcal{F}_h^{-1} & \& & \mathcal{F}_h\circ\left(e^{-i\theta}(-\Delta_h)^{1-\alpha}\partial_t-\Delta_h\right) \circ \mathcal{F}_h^{-1}.
\end{eqnarray*}

Indeed, from (\ref{CliffordWittBasis}) we obtain
\begin{eqnarray*}
	\textbf{z}_{h}(\xi)\f+\f\textbf{z}_{h}(\xi)&=&\textbf{z}_{h}(\xi)\f^\dagger+\f^\dagger\textbf{z}_{h}(\xi)=0\\ \ \\
	\left(\f(d_h(\xi)^2)^{1-\alpha} \partial_t +\f^\dagger e^{-i\theta}\right)^2 &=&\f \f^\dagger e^{-i\theta}(d_h(\xi)^2)^{1-\alpha} \partial_t +\f^\dagger \f e^{-i\theta}(d_h(\xi)^2)^{1-\alpha} \partial_t \\
	&=&e^{-i\theta}(d_h(\xi)^2)^{1-\alpha} \partial_t
\end{eqnarray*}
so that
\begin{eqnarray*}
	\left(\textbf{z}_{h}(\xi)+ \f(d_h(\xi)^2)^{1-\alpha} \partial_t +\f^\dagger e^{-i\theta}\right)^2
	&=&\textbf{z}_{h}(\xi)^2+\left(\f(d_h(\xi)^2)^{1-\alpha} \partial_t +\f^\dagger e^{-i\theta}\right)^2 \\
	&+& \textbf{z}_{h}(\xi)\left(\f(d_h(\xi)^2)^{1-\alpha} \partial_t +\f^\dagger e^{-i\theta}\right)+\left(\f(d_h(\xi)^2)^{1-\alpha} \partial_t +\f^\dagger e^{-i\theta}\right)\textbf{z}_{h}(\xi)\\
	&=&
	e^{-i\theta}(d_h(\xi)^2)^{1-\alpha}\partial_t+d_h(\xi)^2.
\end{eqnarray*}

Regarding the definition of $(-\Delta_h)^\sigma$, we would like to stress the tangible connection between eq.~(\ref{fractionalDeltah}) and the Bochner's definition (cf.~\cite[Section 6]{CiaurriGRTV17}). Essentially, from the fact that the Fourier symbol $e^{-td_h(\xi)^2}$ belongs to $C^\infty\left(\Qh;\BC \otimes \cl_{n,n}\right)$ one can moreover show that $(-\Delta_h)^{\sigma}$ admits, for every of $0<\sigma<1$, a {\it semidiscrete heat semigroup representation} in terms of $\left\{e^{s\Delta_h}\right\}_{s\geq 0}$. Namely, a wise adaptation of the proof of \cite[Lemma 6.5.]{LizamaRoncal18} (see also \cite[Section 2.]{CiaurriGRTV18}):
\begin{eqnarray*}
	\left(d_h(\xi)^2\right)^\sigma= \int_{0}^{\infty} g_{-\sigma}(s)\left(e^{-sd_h(\xi)^2}-1\right)ds, &  g_{-\sigma}(s)=\dfrac{s^{-1-\sigma}}{\Gamma(-\sigma)}~  [\mbox{see eq.}~(\ref{GelfandShilov})]
\end{eqnarray*}
allows us to guarantee that $(d_h(\xi)^2)^\sigma$ also belongs to $C^\infty\left(\Qh;\BC \otimes \cl_{n,n}\right)$  
and moreover that $(-\Delta_h)^\sigma$, represented as follows:
\begin{eqnarray*} (-\Delta_h)^\sigma =\left\{\begin{array}{lll} 
		\displaystyle \int_{0}^{\infty} g_{-\sigma}(s)\left(e^{s\Delta_h}-I\right)ds&,~~0<\sigma<1
		\\ \ \\
		-\Delta_h &,~~\sigma=1
	\end{array}\right.
\end{eqnarray*}
satisfies the mapping property
$$(-\Delta_h)^{\sigma}:\mathcal{S}(h\BZ^n;\BC \otimes \cl_{n,n})\rightarrow \mathcal{S}(h\BZ^n;\BC \otimes \cl_{n,n}).$$

Thus, the space-fractional regularization ${~}_\theta\mathcal{D}_{h,t}^\alpha$ of ${~}_\theta\mathcal{D}_{h,t}$, defined viz eq. (\ref{SpaceFractionalDh}), is thus well defined.
The foregoing lemma permits us to derive the following integro-differential-difference representations, involving the {\it space-fractional semidiscrete} operators ${~}_\theta\mathcal{D}_{h,t}^\alpha$ and $e^{-i\theta}\partial_t(-\Delta_h)^{1-\alpha}-\Delta_h$, respectively.
\begin{lemma}
	For every $t\in [0,\infty)$, let us assume that the components $\Psi^{[m]}(x,t)$ of the
	$\BC\otimes\cl_{n+1,n+1}-$valued function  \begin{eqnarray*}
		\Psi(x,t)=\Psi^{[0]}(x,t)+\f\Psi^{[1]}(x,t)+\f^\dagger\Psi^{[2]}(x,t)+\f\f^\dagger\Psi^{[3]}(x,t)
	\end{eqnarray*}
	satisfy the set of conditions ($m=0,1,2,3$)
	\begin{eqnarray*}
		\Psi^{[m]}(\cdot,t)\in \mathcal{S}(h\BZ^n;\BC \otimes \cl_{n,n}) &\&& \partial_t\Psi^{[m]}(\cdot,t)\in \mathcal{S}(h\BZ^n;\BC \otimes \cl_{n,n}). 
	\end{eqnarray*}

	Then, for every $0<\alpha\leq 1$ and $|\theta|\leq \frac{\alpha \pi}{2}$, the componentwise action
	$${~}_\theta\mathcal{D}_{h,t}^\alpha\Psi(x,t)=D_h\Psi(x,t)+ \f~(-\Delta_h)^{1-\alpha}\partial_t \Psi(x,t)+\f^\dagger e^{-i\theta}\Psi(x,t)$$ admits the following integro-differential-difference representation
	\begin{eqnarray*}
		\label{SpaceFractionalDirac} \left\{\begin{array}{lll} 
			D_h\Psi(x,t)+\displaystyle \int_{0}^{\infty} g_{\alpha-1}(s)\left(e^{s\Delta_h}-I\right)\f\partial_t\Psi(x,t)ds+\f^\dagger e^{-i\theta}\Psi(x,t),~~0<\alpha<1 &
			\\ \ \\
			D_h\Psi(x,t)+\f\partial_t\Psi(x,t)+\f^\dagger e^{-i\theta}\partial_t\Psi(x,t) ,~~\alpha=1,
		\end{array}\right.
	\end{eqnarray*}
	whereas
	\begin{eqnarray*}
		\label{SpaceFractionalLaplacian} \left\{\begin{array}{lll} 
			\displaystyle \int_{0}^{\infty} g_{\alpha-1}(s)\left(e^{s\Delta_h}-I\right)e^{-i\theta}\partial_t\Psi(x,t)ds-\Delta_h\Psi(x,t)&,~~0<\alpha<1
			\\ \ \\
			e^{-i\theta}\partial_t\Psi(x,t)-\Delta_h\Psi(x,t) &,~~\alpha=1
		\end{array}\right.
	\end{eqnarray*}
	stands for the integro-differential-difference representation of $$e^{-i\theta}(-\Delta_h)^{1-\alpha}\partial_t\Psi(x,t)-\Delta_h\Psi(x,t)={~}_\theta\mathcal{D}_{h,t}^\alpha\left({~}_\theta\mathcal{D}_{h,t}^\alpha\Psi(x,t)\right).$$
\end{lemma}

Similarly to {\bf Theorem \ref{CoupledSystemParabolicDiracRLt}}, one can also obtain the following characterization for the {\it space-fractional semidiscrete operator} 
${~}_\theta\mathcal{D}_{h,t}^\alpha$ defined through eq. eq.~(\ref{SpaceFractionalDh}). Although the definition of ${~}_\theta\mathcal{D}_{h,t}^\alpha$ relies essentially on the replacement $$-e^{i \theta(1-\beta)}~\RL^\beta\longrightarrow (-\Delta_h)^{1-\alpha}\partial_t$$
on the right-hand side of (\ref{TimeFractionalDh}), its proof involves a lot of technicalities far beyond the Schwarz class $\mathcal{S}(h\BZ^n;\BC \otimes \cl_{n,n})$.
\begin{theorem}\label{CoupledSystemParabolicDiract}
	For every $t\in [0,\infty)$, let us assume that the components $\Psi^{[m]}(x,t)$ of the
	$\BC\otimes\cl_{n+1,n+1}-$valued function \begin{eqnarray*}
		\Psi(x,t)=\Psi^{[0]}(x,t)+\f\Psi^{[1]}(x,t)+\f^\dagger\Psi^{[2]}(x,t)+\f\f^\dagger\Psi^{[3]}(x,t),
	\end{eqnarray*}
	satisfy	the set of conditions ($m=0,1,2,3$)
	\begin{eqnarray*}
		\Psi^{[m]}(\cdot,t)\in \mathcal{S}(h\BZ^n;\BC \otimes \cl_{n,n}) &\&& \partial_t\Psi^{[m]}(\cdot,t)\in \mathcal{S}(h\BZ^n;\BC \otimes \cl_{n,n}). 
	\end{eqnarray*}

	Then, for every $0<\alpha\leq 1$ and $|\theta|\leq\frac{\alpha\pi}{2}$ the function $\Psi(x,t)$ is a null solution of ${~}_\theta\mathcal{D}_{h,t}^{\alpha}$ (see eq.~(\ref{SpaceFractionalDh})) if, and only if
	\begin{eqnarray}
		\label{CoupledSystemsHeatAlpha} \left\{\begin{array}{lll} 
			\partial_t\Psi^{[m]}(x,t)=-e^{i\theta}(-\Delta_h)^\alpha\Psi^{[m]}(x,t) & , & m=0,2
			\\ \ \\
			\Psi^{[1]}(x,t)=-e^{i\theta}D_h
			\Psi^{[0]}(x,t) & &\\ \ \\
			\Psi^{[3]}(x,t)=e^{i\theta}D_h\Psi^{[2]}(x,t)-\Psi^{[0]}(x,t). &  & 
		\end{array}\right.
	\end{eqnarray}
\end{theorem}

\proof
First, we recall that from (\ref{DiracEqh}) and (\ref{WittBasis}) it readily follows that the discrete Dirac operator $D_h$ and the Witt basis $\f,\f^\dagger$ of $\cl_{1,1}$ satisfy the set of relations
\begin{eqnarray*}
	\f ~D_h=-D_h~\f, & \f^\dagger ~D_h=-D_h~\f^\dagger &~~\mbox{\&}~~ ~D_h~\f\f^\dagger=\f\f^\dagger~ D_h.
\end{eqnarray*}

Then, by letting act the operator (\ref{SpaceFractionalDh}) on $\Psi(x,t)$ one gets, by a straightforwardly computation based on the aforementioned relations, that
\begin{eqnarray*}
	{~}_\theta \mathcal{D}_{h,t}^\alpha\Psi(x,t)&=&	\left(D_h+\f (-\Delta_h)^{1-\alpha}\partial_t+ \f^\dagger~e^{-i\theta}\right)\Psi(x,t) \\
	&=&D_h\Psi^{[0]}(x,t)-\f D_h\Psi^{[1]}(x,t)\\
	&-&\f^\dagger D_h\Psi^{[2]}(x,t)+\f\f^\dagger D_h\Psi^{[3]}(x,t) 
	\\
	&+&\f(-\Delta_h)^{1-\alpha}\partial_t\Psi^{[0]}(x,t)\\
	&+&\f\f^\dagger(-\Delta_h)^{1-\alpha}\partial_t\Psi^{[2]}(x,t)\\
	&+&\f^\dagger\left(e^{-i\theta}\Psi^{[0]}(x,t)+e^{-i\theta}\Psi^{[3]}(x,t)\right)\\
	&+&(1-\f\f^\dagger)e^{-i\theta}\Psi^{[1]}(x,t).
\end{eqnarray*}

By rearranging now the previous identity we obtain the above ansatz, written as a linear combination in terms of $1,\f,\f^\dagger$ and $\f\f^\dagger$. Namely, one has
\begin{eqnarray*}
	{~}_\theta \mathcal{D}_{h,t}^\alpha\Psi(x,t) = \left(D_h\Psi^{[0]}(x,t)+e^{-i\theta}\Psi^{[1]}(x,t)\right) \\
	+\f\left((-\Delta_h)^{1-\alpha}\partial_t\Psi^{[0]}(x,t)-D_h\Psi^{[1]}(x,t)\right)\\
	+\f^\dagger\left(-D_h\Psi^{[2]}(x,t)+e^{-i\theta}\Psi^{[0]}(x,t)+e^{-i\theta}\Psi^{[3]}(x,t)\right)
	\\
	+\f\f^\dagger\left(\partial_t(-\Delta_h)^{1-\alpha}\Psi^{[2]}(x,t)+D_h\Psi^{[3]}(x,t)-e^{-i\theta}\Psi^{[1]}(x,t)\right).
\end{eqnarray*}

Then, simply the observation that $$(-\Delta_h)^{\alpha-1}=\mathcal{F}_h^{-1}\circ (d_h(\xi)^2)^{\alpha-1}\circ \mathcal{F}_h$$ stands for the inverse of $(-\Delta_h)^{1-\alpha}=\mathcal{F}_h^{-1}\circ (d_h(\xi)^2)^{1-\alpha}\circ \mathcal{F}_h$ (see eq. (\ref{fractionalDeltah})), one can thus show that $\Psi(x,t)$ solves the equation ${~}_\theta \mathcal{D}_{h,t}^\alpha\Psi(x,t)=0$ if, and only if
\begin{eqnarray}
	\label{CoupledSystemsDiracAlpha} \left\{\begin{array}{lll} 
		\Psi^{[1]}(x,t)=-e^{i\theta}D_h\Psi^{[0]}(x,t) & (\mbox{eq.}~1) &
		\\ \ \\
		\partial_t\Psi^{[0]}(x,t)=(-\Delta_h)^{\alpha-1}D_h\Psi^{[1]}(x,t) & (\mbox{eq.}~\f)&\\ \ \\
		\Psi^{[3]}(x,t)=e^{i\theta}D_h\Psi^{[2]}(x,t)-\Psi^{[0]}(x,t) &   (\mbox{eq.}~\f^\dagger) & \\ \ \\
		\partial_t\Psi^{[2]}(x,t)=-(-\Delta_h)^{\alpha-1}D_h\Psi^{[3]}(x,t)+e^{-i\theta}(-\Delta_h)^{\alpha-1}\Psi^{[1]}(x,t) &   (\mbox{eq.}~\f\f^\dagger). &
	\end{array}\right.
\end{eqnarray}

Next, using the fact that the discrete Dirac operator $D_h$ (see eq. (\ref{DiracEqh})) satisfies $(D_h)^2=-\Delta_h$ (cf.~\cite[Proposition 2.1]{FaustinoKGordonDirac16}), we immediately get
\begin{eqnarray*}
	\partial_t\Psi^{[0]}(x,t)&=&-e^{i\theta}(-\Delta_h)^{\alpha-1}(D_h)^2\Psi^{[0]}(x,t) \\
	&=&-e^{i\theta}(-\Delta_h)^\alpha\Psi^{[0]}(x,t),
\end{eqnarray*}
after susbtituting
$(\mbox{eq.}~1)$ on the right-hand side of $(\mbox{eq.}~\f)$.

In the same order of ideas, by substituting $(\mbox{eq.}~1)$ \& $(\mbox{eq.}~\f^\dagger)$ on the right-hand side of $(\mbox{eq.}~\f\f^\dagger)$, we end up with
\begin{eqnarray*}
	\partial_t\Psi^{[2]}(x,t)&=&-(-\Delta_h)^{\alpha-1}D_h\left(e^{i\theta}D_h\Psi^{[2]}(x,t)-\Psi^{[0]}(x,t)\right)+e^{-i\theta}(-\Delta_h)^{\alpha-1}\left(-e^{i\theta}D_h\Psi^{[0]}(x,t)\right) \\
	&=&-e^{i\theta}(-\Delta_h)^{\alpha-1}(D_h)^2\Psi^{[2]}(x,t) \\
	&=&-e^{i\theta}(-\Delta_h)^{\alpha}\Psi^{[2]}(x,t).
\end{eqnarray*}

Thus, we have shown that the coupled system of equations (\ref{CoupledSystemsDiracAlpha}) is equivalent to (\ref{CoupledSystemsHeatAlpha}).
\qed

\begin{remark}
	Under the conditions of {\bf Theorem \ref{CoupledSystemParabolicDiract}}, one can say that the $\BC\otimes\cl_{n,n}-$valued components $\Psi^{[0]}(x,t)$ and $\Psi^{[2]}(x,t)$ of $\Psi(x,t)$ are null solutions of the {\it space-fractional semidiscrete operator} $e^{-i\theta}\partial_t+(-\Delta_h)^\alpha$.
\end{remark}

\section{Main Results}\label{MainResultsSection}

\subsection{The analytic semigroup $\left\{\exp\left(-te^{i\theta}(-\Delta_h)^\alpha\right)\right\}_{t\geq 0}$}
In this section we will show that the action of the [fractional semidiscrete] analytic semigroup $\left\{\exp\left(-te^{i\theta}(-\Delta_h)^\alpha\right)\right\}_{t\geq 0},$ carrying the parameters $0<\alpha\leq 1$ and $|\theta|\leq \frac{\alpha \pi}{2}$, allows us to establish a correspondence between the solution representation of two seemingly distinct Cauchy problems, involving the space-fractional and time-fractional semidiscrete Dirac operators studied previously in Section \ref{SemidiscreteSection}.

Before we proceed, it is important to recall some facts regarding the formulation of $\left\{\exp(-te^{i\theta}(-\Delta_h)^\alpha)\right\}_{t\geq 0}$ as an analytic semigroup encoded by the product of {\it modified Bessel functions of the first kind} along the same lines of \cite[Subsection 21.4.2]{F19}.
First, notice that in case of $\alpha=1$, $\left\{\exp\left(te^{i\theta}\Delta_h\right)\right\}_{t\geq 0}$ corresponds to an analytic extension of the {\it semidiscrete heat semigroup}, already treated in \cite[Section 6]{LizamaRoncal18} (see {\bf Remark 13}). In concrete, from the \textit{discrete 
	convolution formula} (\ref{discreteConvolution})
there holds
\begin{eqnarray}
	\label{discreteConvolutionHeat}\exp\left(te^{i\theta}\Delta_h\right)\Phi(x)&=&\sum_{y \in h\BZ^n}  h^n\Phi(y)K(x-y,te^{i\theta}),
\end{eqnarray}
with 
\begin{eqnarray}
	\label{IntegralDiscreteHeat}K(x,te^{i\theta})=\frac{1}{(2\pi)^{n}} \int_{Q_h} e^{-te^{i\theta}d_h(\xi)^2} e^{-i x \cdot \xi} d\xi.
\end{eqnarray}
Moreover, the closed formula\footnote{The published version of the book chapter \cite{F19} contains a constant typo in the formula (\ref{IntegralDiscreteHeat}), which has been corrected here.} (cf.~\cite[p. 458, eq.~(21.35)]{F19})
\begin{eqnarray*}
	\label{ClosedFormulaBessel} K(x,te^{i\theta})&=&\frac{1}{h^{n}}e^{-\frac{2nte^{i\theta}}{h^2}} I_{\frac{x_1}{h}}\left(\frac{2te^{i\theta}}{h^2}\right)I_{\frac{x_2}{h}}\left(\frac{2te^{i\theta}}{h^2}\right) \ldots I_{\frac{x_n}{h}}\left(\frac{2te^{i\theta}}{h^2}\right),
\end{eqnarray*}
involving the product of {\it modified Bessel functions of the first kind} (cf.~\cite[p. 456, 2.5.40 (3)]{PBM86})
\begin{eqnarray}
	\label{ModifiedBessel}I_k(z)=\frac{1}{\pi}\int_{0}^\pi e^{z\cos(\omega)}e^{-ik\omega}d\omega, & |\mbox{arg}(z)|<\pi.
\end{eqnarray}
results from the identity associated to the Fourier multipliers $d_h(\xi)^2$ (see eq.~(\ref{FourierMultipliers})):
$$ 
d_h(\xi)^2=\sum_{j=1}^n \frac{2}{h^2}\left(1-\cos(h\xi_j)\right),
$$
and from the change of variables $\xi_j=\frac{\omega_j}{h}$ ($-\pi <\omega_j\leq \pi$) on (\ref{IntegralDiscreteHeat}).

Because of $|\mbox{arg}(z)|=|\theta|$ for $z=te^{i\theta}$ one can moreover say that $K(x,te^{i\theta})$, described as above, is well defined for every $|\theta|\leq \frac{\alpha \pi}{2}$ ($\alpha=1$). This conclusion is then immediate from the integral representation~(\ref{ModifiedBessel}).

We can moreover extend the above analysis to $\left\{\exp\left(-te^{i\theta}(-\Delta_h)^\alpha\right)\right\}_{t\geq 0}$ in case of $0<\alpha<1$ and $|\theta|\leq \frac{\alpha \pi}{2}$. At this stage, we would like to stress that that the Laplace identity
\begin{eqnarray}
	\label{Levy}e^{-s^\alpha}=\int_{0}^\infty e^{-rs}L_\alpha^{-\alpha}(r)dr, & \Re(s)>0, & 0<\alpha<1
\end{eqnarray}
involving the {\it L\'evy stable distribution} $L_\alpha^{-\alpha}(r)$ (cf.~\cite[Section 4]{Mainardi01}) provides us to establish a bridge between the semigroups
\begin{center}
	$\left\{\exp\left(-te^{i\theta}(-\Delta_h)^\alpha\right)\right\}_{t\geq 0}$ and $\left\{\exp\left(te^{\frac{i\theta}{\alpha}}\Delta_h\right)\right\}_{t\geq 0}$,
\end{center}
respectively.

Indeed, for $s=te^{i\theta}(d_h(\xi)^2)^\alpha$ one can see from the change of variable $r\rightarrow rt^{-\frac{1}{\alpha}}$ on the right-hand side of (\ref{Levy}) that the underlying Fourier symbols $e^{-te^{i\theta}(d_h(\xi)^2)^\alpha}$ and $e^{-te^{\frac{i\theta}{\alpha}}d_h(\xi)^2}$, respectively, are linked by the integral formula
\begin{eqnarray*}
	e^{-te^{i\theta}(d_h(\xi)^2)^\alpha}=\int_{0}^{\infty} e^{-re^{\frac{i\theta}{\alpha}}d_h(\xi)^2}f_{\alpha,\theta}(r) du,~~f_{\alpha,\theta}(r)=t^{-\frac{1}{\alpha}}L_\alpha^{-\alpha}\left(rt^{-\frac{1}{\alpha}}\right),
\end{eqnarray*}
whereby $t>0$ is assumed to ensure that the constraint $\mbox{Re}(s)>0$, appearing on eq.~(\ref{Levy}), is always fulfilled.

The above representation is useful in several applications on the crossroads of fractional calculus and stochastics, but we will not need to apply it in concrete throughout the present paper. For the interested reader, we refer e.g. to \cite[Chapter 3 \& Chapter 6]{MS12}.

\subsection{Space-fractional vs time-fractional Cauchy problems}

For the remainder part of this paper we will restrict ourselves to the analysis of $\left\{\exp\left(-te^{i\theta}(-\Delta_h)^\alpha\right)\right\}_{t\geq 0}$ in terms of its symbol $e^{-te^{i\theta}(d_h(\xi)^2)^\alpha}$. 
First of all, we will start to show that the technique used in \cite{CLRV15} to prove {\bf Theorem 3} can be generalized to the fractional semidiscrete analytic semigroup $\left\{\exp\left(-te^{i\theta}(-\Delta_h)^\alpha\right)\right\}_{t\geq 0}$.
\begin{theorem}\label{SemigroupCauchySolver}
	Let $\Phi_0\in \mathcal{S}(h\BZ^n;\BC \otimes \cl_{n,n})$ be given. Then, 
	for every $0<\alpha\leq 1$ and $|\theta|\leq \frac{\alpha \pi}{2}$ the function 
	$$
	\Phi(x,t)=\exp(-te^{i\theta}(-\Delta_h)^\alpha)\Phi_0(x)
	$$
	solves the following two Cauchy problems:
	\begin{eqnarray}
		\label{CauchySemidiscrete} \left\{\begin{array}{lll} 
			\partial_t\Phi(x,t)= -e^{i\theta }(-\Delta_h)^\alpha \Phi(x,t) & \mbox{for} & (x,t)\in
			h\BZ^n \times [0,\infty)
			\\ \ \\
			\Phi(x,0)=\Phi_0(x) & \mbox{for} & x\in h\BZ^n,
		\end{array}\right.
	\end{eqnarray}
	and
	\begin{eqnarray}
		\label{CauchyTimeFractional} \left\{\begin{array}{lll} 
			\mathbb{D}_t^{\frac{1}{\alpha}}\Phi(x,t)= -e^{\frac{i\theta}{\alpha}}\Delta_h \Phi(x,t) & \mbox{for} & (x,t)\in
			h\BZ^n \times [0,\infty)
			\\ \ \\
			\Phi(x,0)=\Phi_0(x) & \mbox{for} & x\in h\BZ^n.
		\end{array}\right.
	\end{eqnarray}
\end{theorem}

\proof
By letting act the discrete Fourier transform $\mathcal{F}_h$, we obtain an equivalent formulation for the Cauchy problem (\ref{CauchySemidiscrete}) on the momentum space $Q_h\times [0,\infty)$:
\begin{eqnarray}
	\label{CauchySemidiscreteFourier} \left\{\begin{array}{lll} 
		\partial_t\left[\mathcal{F}_h\Phi(\xi,t)\right]= -e^{i\theta}(d_h(\xi)^2)^{\alpha} \mathcal{F}_h\Phi(\xi,t) & \mbox{for} & (\xi,t)\in
		Q_h \times [0,\infty)
		\\ \ \\
		\mathcal{F}_h\Phi(\xi,0)=\mathcal{F}_h\Phi_0(\xi) & \mbox{for} & \xi \in Q_h
	\end{array}\right.
\end{eqnarray}
so that 
\begin{eqnarray}
	\label{SolutionSemidiscreteFourier}\mathcal{F}_h\Phi(\xi,t)= e^{-te^{i\theta}\left(d_h(\xi)^2\right)^\alpha} \mathcal{F}_h\Phi_0(\xi)& ,(\xi,t)\in Q_h \times [0,\infty)
\end{eqnarray}
solves (\ref{CauchySemidiscreteFourier}), and whence
\begin{eqnarray*}
	\Phi(x,t)&=&\exp(-te^{i\theta}(-\Delta_h)^\alpha)\Phi_0(x)\\
	&=&\displaystyle \frac{1}{(2\pi)^{\frac{n}{2}}}\int_{Q_h} e^{-te^{i\theta}\left(d_h(\xi)^2\right)^\alpha} \mathcal{F}_h\Phi_0(\xi) e^{-ix\cdot \xi} d\xi
\end{eqnarray*}
solves (\ref{CauchySemidiscrete}).

Thus, in order to show that $\Phi(x,t)$ is also a solution of the Cauchy problem (\ref{CauchyTimeFractional}), it suffices to show that (\ref{SolutionSemidiscreteFourier}) solves 
\begin{eqnarray}
	\label{CauchyTimeFractionalFourier} \left\{\begin{array}{lll} 
		\mathbb{D}_t^{\frac{1}{\alpha}}\left[\mathcal{F}_h\Phi(\xi,t)\right]=e^{\frac{i\theta}{\alpha}} d_h(\xi)^2\mathcal{F}_h\Phi(\xi,t) & \mbox{for} & (\xi,t)\in Q_h \times [0,\infty)
		\\ \ \\
		\mathcal{F}_h\Phi(\xi,0)=\mathcal{F}_h\Phi_0(\xi) & \mbox{for} & \xi\in Q_h,
	\end{array}\right.
\end{eqnarray}
or equivalently, that $$\mathbb{D}_t^{\frac{1}{\alpha}}\left[e^{-te^{i\theta}\left(d_h(\xi)^2\right)^\alpha}\right]=e^{\frac{i\theta}{\alpha}}d_h(\xi)^2~e^{-te^{i\theta}\left(d_h(\xi)^2\right)^\alpha},$$
holds for every $0<\alpha\leq 1$ and $|\theta|\leq \frac{\alpha \pi}{2}$.

To do so, recall that for every $k\in \mathbb{N}$ one has the derivation rule
\begin{eqnarray}
	\label{dkExp}\left(-\partial_t\right)^{k}\left[e^{-te^{i\theta}\left(d_h(\xi)^2\right)^\alpha}\right]=e^{i \theta  k}\left(d_h(\xi)^2\right)^{\alpha k}e^{-te^{i\theta}\left(d_h(\xi)^2\right)^\alpha}. 
\end{eqnarray}

If $\frac{1}{\alpha}=k\in \BN$, it readily follows that $\mathbb{D}_t^{\frac{1}{\alpha}}$ equals to $(-\partial_t)^k$ so that 
\begin{eqnarray*}
	\mathbb{D}_t^{\frac{1}{\alpha}}\left[e^{-te^{i\theta}\left(d_h(\xi)^2\right)^\alpha} \right]=e^{\frac{i\theta}{\alpha}}
	d_h(\xi)^2e^{-te^{i\theta}\left(d_h(\xi)^2\right)^\alpha}.
\end{eqnarray*}

Otherwise, set $k=\left\lfloor \frac{1}{\alpha} \right\rfloor+1$. Note that, for values of $0<\alpha<1$, the constraint $|\theta|\leq \frac{\alpha \pi}{2}$ assures that the constant $\lambda=e^{i\theta}\left(d_h(\xi)^2\right)^\alpha$ satisfies the condition
$$
\Re(\lambda)=\cos(\theta)\left(d_h(\xi)^2\right)^\alpha>0.
$$

Then, from the Laplace identity (\ref{LaplaceId}) we get that 
$$
\int_0^\infty g_{k-\frac{1}{\alpha}}(p)e^{-pe^{i\theta}\left(d_h(\xi)^2\right)^\alpha}~dp=e^{-i \theta\left(k-\frac{1}{\alpha}\right)}\left(d_h(\xi)^2\right)^{-\alpha\left(k-\frac{1}{\alpha}\right)},
$$
whereby $g_{k-\frac{1}{\alpha}}(p)$ stands for the Gel'fand-Shilov function (see eq.~(\ref{GelfandShilov})).

Hence, the eq. (\ref{dkExp}) together with the previous integral formula gives rise to the sequence of identities
\begin{eqnarray*}
	\mathbb{D}_t^{\frac{1}{\alpha}}\left[e^{-te^{i\theta }\left(d_h(\xi)^2\right)^\alpha}\right]&=&\left(-\partial_t\right)^{k} \int_t^\infty g_{k-\frac{1}{\alpha}}(s-t)e^{-se^{i\theta}\left(d_h(\xi)^2\right)^\alpha}~ds \\
	&=&\left(-\partial_t\right)^{k} \int_0^\infty g_{k-\frac{1}{\alpha}}(p)~e^{-(t+p)e^{i\theta}\left(d_h(\xi)^2\right)^\alpha}~dp \\
	&=&\left(-\partial_t\right)^{k}\left[e^{-te^{i\theta}\left(d_h(\xi)^2\right)^\alpha}\right]\times  \int_0^\infty g_{k-\frac{1}{\alpha}}(p)~e^{-pe^{i\theta}\left(d_h(\xi)^2\right)^\alpha}~dp\\
	&=&e^{i \theta k}\left(d_h(\xi)^2\right)^{\alpha k}e^{-te^{i\theta}\left(d_h(\xi)^2\right)^\alpha} \times e^{-i \theta\left(k-\frac{1}{\alpha}\right)}\left(d_h(\xi)^2\right)^{-\alpha(k-\frac{1}{\alpha})} \\
	&=&e^{\frac{i \theta}{\alpha}}d_h(\xi)^2~e^{-te^{i\theta}\left(d_h(\xi)^2\right)^\alpha},
\end{eqnarray*}
concluding the proof of (\ref{CauchyTimeFractionalFourier}).
\qed
\begin{remark}
	Essentially, we have shown in {\bf Theorem \ref{SemigroupCauchySolver}} that the analytic semigroup $\left\{\exp\left(-te^{i\theta}(-\Delta_h)^\alpha\right)\right\}_{t\geq 0}$, carrying the parameters $0<\alpha\leq 1$ and $|\theta|\leq\frac{\alpha\pi}{2}$, generates simultaneously solutions for Cauchy problems induced by the fractional semidiscrete operators $e^{-i\theta}\partial_t+(-\Delta_h)^\alpha$ and $-e^{-\frac{i\theta}{\alpha}}\RL^{\frac{1}{\alpha}}-\Delta_h$, respectively.
\end{remark}

\subsection{Cauchy problems of L\'evy-Leblond type}

After the meaningful construction obtained previously in {\bf Theorem \ref{SemigroupCauchySolver}}, we have now gathered the main ingredients to prove the main result of this section, whose starting point relies heavily on 
{\bf Theorem \ref{CoupledSystemParabolicDiracRLt}} \& {\bf Theorem \ref{CoupledSystemParabolicDiract}}.

\begin{theorem}\label{LevyLeblondSolutions}
	Let $\Phi^{[0]}_0,\Phi^{[2]}_0\in \mathcal{S}(h\BZ^n;\BC\otimes \cl_{n,n})$ and set
	\begin{eqnarray*}
		\Phi_0(x)&=&\f^\dagger\f\Phi^{[0]}_0(x)+\f^\dagger\Phi^{[2]}_0(x), \\
		\Phi(x,t)&=&\exp(-te^{i\theta}(-\Delta_h)^\alpha)\Phi_0(x).
	\end{eqnarray*}
	
	Then, for every $0<\alpha\leq 1$ and $|\theta|\leq \frac{\alpha \pi}{2}$ the function
	\begin{eqnarray}
		\label{HarmonicConjugatesTheta}\Psi(x,t)&=&\Phi(x,t)-\f e^{i\theta}D_h\Phi(x,t)
	\end{eqnarray}
	solves the following two Cauchy problems of L\'evy-Leblond type: \begin{eqnarray}
		\label{CauchyFractionalDiracAlpha} \left\{\begin{array}{lll} 
			\left(D_h+\f (-\Delta_h)^{1-\alpha}\partial_t+ \f^\dagger~e^{-i\theta}\right)\Psi(x,t)=0 & \mbox{for} & (x,t)\in
			h\BZ^n \times [0,\infty)
			\\ \ \\
			\Psi(x,0)=\Phi_0(x)
			-\f e^{i\theta}D_h\Phi_0(x)& \mbox{for} & x\in h\BZ^n,
		\end{array}\right.
	\end{eqnarray}
	and
	\begin{eqnarray}
		\label{CauchyFractionalRLt} \left\{\begin{array}{lll} 
			\left(D_h-\f~e^{i\theta \left(1-\frac{1}{\alpha}\right)} \RL^{\frac{1}{\alpha}}+ \f^\dagger~e^{-i\theta}\right)\Psi(x,t)=0 & \mbox{for} & (x,t)\in
			h\BZ^n \times [0,\infty)
			\\ \ \\
			\Psi(x,0)=\Phi_0(x)
			-\f e^{i\theta}D_h\Phi_0(x)& \mbox{for} & x\in h\BZ^n.
		\end{array}\right.
	\end{eqnarray}
	
\end{theorem}

\proof Under the assumption that $\Phi^{[0]}_0,\Phi^{[2]}_0\in \mathcal{S}(h\BZ^n;\BC\otimes \cl_{n,n})$, {\bf Theorem \ref{SemigroupCauchySolver}} asserts that
the components 
\begin{eqnarray*}
	\Psi^{[m]}(x,t)=\exp\left(-te^{i\theta}(-\Delta_h)^\alpha\right)\Phi_0^{[m]}(x) & (m=0,2)
\end{eqnarray*}
of
\begin{eqnarray*}
	\Phi(x,t)=\exp\left(-te^{i\theta}(-\Delta_h)^\alpha\right)\Phi_0(x), &\mbox{with} & \Phi_0(x)=\f^\dagger\f\Phi^{[0]}_0(x)+\f^\dagger\Phi^{[2]}_0(x)
\end{eqnarray*}
are solutions of the following set of Cauchy problems ($m=0,2$) for values of $0<\alpha\leq 1$ and $|\theta|\leq \frac{\alpha \pi}{2}$:
\begin{eqnarray}
	\label{CauchySemidiscreteComponentsAlpha} \left\{\begin{array}{lll} 
		\partial_t\Psi^{[m]}(x,t)= -e^{i\theta}(-\Delta_h)^\alpha \Psi^{[m]}(x,t) & \mbox{for} & (x,t)\in
		h\BZ^n \times [0,\infty)
		\\ \ \\
		\Psi^{[m]}_0(x,0)=\Phi_0^{[m]}(x) & \mbox{for} & x\in h\BZ^n,
	\end{array}\right.
\end{eqnarray}
and \begin{eqnarray}
	\label{CauchySemidiscreteComponentsAlphaRL} \left\{\begin{array}{lll} 
		\mathbb{D}_t^{\frac{1}{\alpha}}\Psi^{[m]}(x,t)= -e^{ \frac{i\theta}{\alpha}}\Delta_h \Psi^{[m]}(x,t) & \mbox{for} & (x,t)\in
		h\BZ^n \times [0,\infty)
		\\ \ \\
		\Psi^{[m]}_0(x,0)=\Phi_0^{[m]}(x) & \mbox{for} & x\in h\BZ^n.
	\end{array}\right.
\end{eqnarray}

Then, from {\bf Theorem \ref{CoupledSystemParabolicDiract}} and {\bf Theorem \ref{CoupledSystemParabolicDiracRLt}} one has that the function 
\begin{eqnarray*}
	\Psi(x,t)=\Psi^{[0]}(x,t)+\f\Psi^{[1]}(x,t)+\f^\dagger\Psi^{[2]}(x,t)+\f\f^\dagger\Psi^{[3]}(x,t),
\end{eqnarray*}
solve simultaneously the Cauchy problems
(\ref{CauchyFractionalDiracAlpha}) and (\ref{CauchyFractionalRLt}) if, and only if, the components $\Psi^{[m]}(x,t)$ ($m=1,3$) of $\Psi(x,t)$ are uniquely determined by 
\begin{eqnarray*}
	\Psi^{[1]}(x,t)&=&-e^{i\theta}D_h\Psi^{[0]}(x,t), \\ \Psi^{[3]}(x,t)&=&e^{i\theta}D_h\Psi^{[2]}(x,t)-\Psi^{[0]}(x,t).
\end{eqnarray*}

Hence, from 
the set of properties 
\begin{eqnarray*}
	\f ~D_h=-D_h~\f, & \f^\dagger ~D_h=-D_h~\f^\dagger &~~\mbox{\&}~~ ~D_h~\f\f^\dagger=\f\f^\dagger~ D_h
\end{eqnarray*}
that yield from the combination of (\ref{DiracEqh}) and (\ref{WittBasis}), it readily follows that
\begin{eqnarray*}
	\Psi(x,t)&=& \Psi^{[0]}(x,t)+\f^\dagger \Psi^{[2]}(x,t) \\
	&+& \f (-e^{i\theta}D_h\Psi^{[0]}(x,t))+\f \f^\dagger\left(e^{i\theta}D_h\Psi^{[2]}(x,t)-\Psi^{[0]}(x,t)\right)\\
	&=&(1-\f \f^\dagger)\Psi^{[0]}(x,t)+\f^\dagger\Psi^{[2]}(x,t)\\
	&+& e^{i\theta}D_h\left(\f \Psi^{[0]}(x,t)\right)+e^{i\theta}D_h\left(\f \f^\dagger  \Psi^{[2]}(x,t)\right)\\
	&=&\f^\dagger \f \Psi^{[0]}(x,t)+\f^\dagger\Psi^{[2]}(x,t)+e^{i\theta}D_h\left[\f\Psi^{[0]}(x,t)+\f\f^\dagger\Psi^{[2]}(x,t)\right],
\end{eqnarray*}
whereas from the identity $\f\f^\dagger \f=\f$, there holds
\begin{eqnarray*}
	\f \Phi(x,t)= \f \left(\f^\dagger \f \Psi^{[0]}(x,t)+ \f^\dagger \Psi^{[2]}(x,t)\right)=\f \Psi^{[0]}(x,t)+\f \f^\dagger \Psi^{[2]}(x,t), \\ \ \\
	-\f e^{i\theta}D_h\Phi(x,t)=e^{i\theta}D_h(\f\Phi(x,t))=e^{i\theta}D_h\left[\f\Psi^{[0]}(x,t)+\f\f^\dagger\Psi^{[2]}(x,t)\right].
\end{eqnarray*}

That allows us to conclude that $\Psi(x,t)$, described as above, is equivalent to (\ref{HarmonicConjugatesTheta}), as desired.
\qed

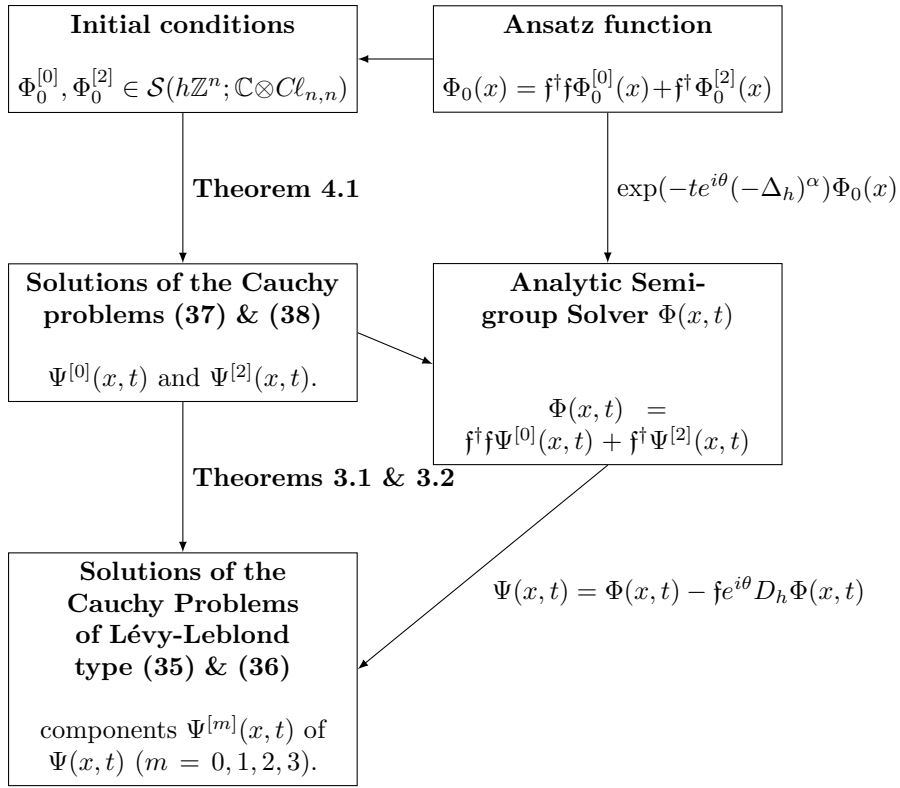
\begin{figure}\label{MainResultsPicture}
	\begin{tikzpicture}
		\node[mynode] (v1){\textbf{Initial conditions}\\ \ \\	$\Phi^{[0]}_0,\Phi^{[2]}_0\in \mathcal{S}(h\BZ^n;\BC\otimes \cl_{n,n})$ };

		\node[mynode, right=1.0cm of v1](v5) {{\bf Ansatz function} \\ \ \\
			$\Phi_0(x)=\f^\dagger \f \Phi_0^{[0]}(x)+\f^\dagger  \Phi_0^{[2]}(x)$	};	
		\node[mynode,below=2.0cm of v5](v3) {{\bf Analytic Semigroup Solver $\Phi(x,t)$} \\ \ \\ \ \\
			$\Phi(x,t)=\f^\dagger \f \Psi^{[0]}(x,t)+\f^\dagger\Psi^{[2]}(x,t)$};
		\node[mynode,below =2.0cm of v1](v2) {
			\textbf{Solutions of the Cauchy problems (\ref{CauchySemidiscreteComponentsAlpha}) \& (\ref{CauchySemidiscreteComponentsAlphaRL})} \\ \ \\ $\Psi^{[0]}(x,t)$ and $\Psi^{[2]}(x,t)$.};
		
		\node[mynode,below = 2.0cm of v2] (v4){{\bf Solutions of the Cauchy Problems of L\'evy-Leblond type (\ref{CauchyFractionalDiracAlpha}) \& (\ref{CauchyFractionalRLt})}\\ \ \\  components $\Psi^{[m]}(x,t)$ of $\Psi(x,t)$~$(m=0,1,2,3)$.};
		
		\draw[-latex] (v5.south) -- node[auto,] {$\exp(-te^{i\theta}(-\Delta_h)^\alpha)\Phi_0(x)$} (v3.north);
		\draw[-latex] (v2.east) -- node[auto,] {} (v3.west);
		\draw[-latex] (v2.south) -- node[auto,] {{\bf Theorems \ref{CoupledSystemParabolicDiracRLt} \& \ref{CoupledSystemParabolicDiract}}} (v4.north);
		\draw[-latex] (v3.south) -- node[auto,] {$\Psi(x,t)=\Phi(x,t)-\f e^{i\theta} D_h\Phi(x,t)$} (v4.east);
		\draw[-latex] (v5.west) -- node[auto,]{} (v1.east);
		\draw[-latex] (v1.south) -- node[auto,]{{\bf Theorem \ref{SemigroupCauchySolver}}} (v2.north);
	\end{tikzpicture}
	\caption{Schematic proof of {\bf Theorem \ref{LevyLeblondSolutions}}.}
\end{figure}

\section{Postscripts}\label{SpaceFractionalRemarks}


\subsection{Factorization of space-fractional semidiscrete operators}
With the proof of {\bf Theorem \ref{LevyLeblondSolutions}}, neatly summarized on \textbf{Figure 2}, it was established an intriguing correspondence between the null solutions of time-fractional resp. space-fractional analogues of the semidiscrete Dirac operator, treated on Subsection \ref{TimeFractionalSection} \& Subsection \ref{SpaceFractionalSection}, and the null solutions of the time-fractional resp. space-fractional regularizations of the semidiscrete Dirac operator (\ref{SemidiscreteDh}).

At this stage we would like to emphasize that, contrary to the time-fractional regularization (\ref{TimeFractionalDh}), the corresponding space-fractional regularization (\ref{SpaceFractionalDh}) of (\ref{SemidiscreteDh}) does not factorize the space-fractional semidiscrete operator $e^{-i\theta}\partial_t+(-\Delta_h)^{\alpha}$, although the factorization of ${~}_\theta\mathbb{D}_{h,t}^{\frac{1}{\alpha}}$ ($\beta=\frac{1}{\alpha}$) and ${~}_\theta\mathcal{D}_{h,t}^{\alpha}$ was never required on the proof of main results of this paper.

To circumvent this gap one can consider the alternative space-fractional semidiscrete variant 
\begin{eqnarray}
	\label{SpaceFractionalVariant}{~}_\theta{\bf D}_{h,t}^\alpha:=(-\Delta_h)^{\frac{\alpha-1}{2}}D_h+\f (-\Delta_h)^{\frac{1-\alpha}{2}}\partial_t+\f^\dagger e^{-i\theta}(-\Delta_h)^{\frac{\alpha-1}{2}}
\end{eqnarray}
of ${~}_\theta\mathcal{D}_{h,t}=D_h+\f\partial_t+\f^\dagger e^{-i\theta}$, involving the space-fractional operator $(-\Delta_h)^{\sigma}$ and its inverse $(-\Delta_h)^{-\sigma}$ ($\sigma=\frac{1-\alpha}{2}$). 

Hereby, the discrete fractional operator $(-\Delta_h)^{-\sigma}D_h$ stands for the hypercomplex extension of the fractional Riesz type transform considered in \cite[Section 6]{CiaurriGRTV17} (see also \cite[Subsection 21.4.3]{F19}).

To do so, let us now turn our attention to the fractional integral operator $(-\Delta_h)^{-\sigma}$. Owing the fact that the Fourier multiplier $(d_h(\xi)^2)^{-\sigma}$ of $\mathcal{F}_h\circ (-\Delta_h)^{-\sigma} \circ \mathcal{F}_h^{-1}$ admitting, for values of $0<\sigma<\frac{1}{2}$,
the following Eulerian integral representation
\begin{eqnarray*}
	(d_h(\xi)^2)^{-\sigma}=\int_0^\infty e^{-pd_h(\xi)^2} g_\sigma(p) dp, & g_\sigma(p)=\dfrac{p^{\sigma-1}}{\Gamma(\sigma)} & [\mbox{see eq.}~(\ref{GelfandShilov})]
\end{eqnarray*}
one can guarantee that the mapping property $$(-\Delta_h)^{-\sigma}:\mathcal{S}(h\BZ^n;\BC \otimes \cl_{n,n})\rightarrow \mathcal{S}(h\BZ^n;\BC \otimes \cl_{n,n}),$$
carrying the inverse of $(-\Delta_h)^{\sigma}$, is fulfilled for every $0<\sigma<\frac{1}{2}$.

As a consequence, the reformulation of $(-\Delta_h)^{\frac{1-\alpha}{2}}$ and $(-\Delta_h)^{\frac{\alpha-1}{2}}$ in terms of the discrete convolution property (\ref{discreteConvolution}):
\begin{eqnarray*}
	(-\Delta_h)^{\sigma}\Psi(x,t) &=& (\Psi (\cdot,t) \star_{h} \delta_{h,\sigma})(x)\\
	&:=&\sum_{y\in h\BZ^n} h^n \delta_{h,\sigma}(x-y)\Psi(y,t),
\end{eqnarray*}
with
$$
\delta_{h,\sigma}(x-y)=
\frac{1}{(2\pi)^{n}}\int_{\Qh} \left(d_h(\xi)^2\right)^\sigma e^{-i(x-y)\cdot \xi}d\xi,
$$
is well defined for every $0<\sigma<\frac{1}{2}$. In case of $\sigma\rightarrow 0^+$ one has that $\delta_{h,\sigma}(x-y)$ converges to the so-called {\it discrete delta function} on $h\BZ^n$. That is,
\begin{eqnarray*}
	\lim_{\sigma\rightarrow 0^+}\delta_{h,\sigma}(x-y)= \left\{\begin{array}{lll} 
		\dfrac{1}{h^n}, &\mbox{if} & x=y
		\\ \ \\
		0 &\mbox{if} & x\neq y.
	\end{array}\right.
\end{eqnarray*}

From the spectral representation $\mathcal{F}_h\circ {~}_\theta{\bf D}_{h,t}^{\alpha}\circ \mathcal{F}_h^{-1}$ in terms of its multiplier,
there holds that
the operators 
(\ref{SpaceFractionalVariant}) and ${~}_\theta\mathcal{ D}_{h,t}^{\alpha}$ (see eq.~(\ref{SpaceFractionalDh})) are interrelated by the formula
$$
{~}_\theta{\bf D}_{h,t}^\alpha=(-\Delta_h)^{\frac{\alpha-1}{2}}{~}_\theta\mathcal{D}_{h,t}^\alpha.
$$

Noteworthy, we obtain the factorization property
$$({~}_\theta{\bf D}_{h,t}^\alpha)^2=e^{-i\theta}\partial_t+(-\Delta_h)^\alpha.$$

\subsection{Function spaces}

In this paper it was provided a successful strategy to solve Cauchy problems envolving discrete space-fractional resp. time-fractional variants of the L\'evy-Leblond operator, also quoted in the literature as {\it parabolic Dirac type operators}. Within this avenue of thought, already considered in the series of papers \cite{F19b,F19,F20}, we have built a discrete pseudo-differential calculus framework in a two fold-way: 
\begin{enumerate}
	\item To embody the multivector structure of Clifford algebras, ubiquitous e.g. on the formulation of discrete boundary value problems of Navier-Stokes (cf.~\cite{faustino2006difference}) and Schr\"odinger type (cf.\cite{CFV08});
	\item To provide a reformulation of the discrete operator calculus considered on the seminal monograph \cite{gurlebeck1997quaternionic} of G\"urlebeck and Spr\"o\ss ig (see \cite[Chapter 5]{gurlebeck1997quaternionic}).
\end{enumerate}

In descriptive terms, the central issue of this paper was to develop a whole machinery in terms of the {\it theory of discrete distributions} over the lattice $h\BZ^n$, bearing in mind that the Schwartz class $\mathcal{S}(h\BZ^n;\BC\otimes \cl_{n,n})$ is dense in $\ell_2(h\BZ^n;\BC\otimes\cl_{n,n})$ (cf.~\cite[Exercise 3.1.15 of p.~302]{RuzhanskyT10}). 

The natural question that naturally arises is the following: {\it 'How the discrete $\ell_p-$spaces, Sobolev spaces and alike, already considered on the monograph \cite{gurlebeck1997quaternionic} and on the papers \cite{faustino2006difference,CFV08}, come into play in this framework?'} 

Firstly, we would like to stress that in stark constrast with the {\it continuum} setting over the Euclidean space $\BR^n$, the {\it discrete Fourier transform} (\ref{discreteFh}) maps the discrete Schwartz space $\mathcal{S}(h\BZ^n;\BC\otimes \cl_{n,n})$ onto the continuous space $C^\infty\left(\Qh;\BC \otimes \cl_{n,n}\right)$. And such isometric isomorphism can be exploited to $L_p-$type spaces under slightly different circumstances, how we can tacitly infer e.g. from the standard proof of Hausdorff–Young inequality (cf.~\cite[Corollary 3.1.24]{RuzhanskyT10}), but also from the proof of embedding result obtained in \cite[LEMMA 3.1]{faustino2006difference} (see also \cite[p.~476]{hytonen2016analysis}).

To be more precise, based on the $L_p-$extension problem posed in \cite[Subsection 2.1.]{hytonen2016analysis} from a probabilistic perspective, one can conjecture that the mapping property $$\mathcal{F}_{h}:\ell_q(h\BZ^n;\BC \otimes \cl_{n,n})\rightarrow L_p\left(\Qh;\BC \otimes \cl_{n,n}\right),$$
yields an isometric isomorphism between the $\BC\otimes\cl_{n,n}-${\it Banach modules} $\displaystyle \ell_q(h\BZ^n;\BC\otimes\cl_{n,n}):=\ell_q(h\BZ^n)\otimes\left(\BC\otimes \cl_{n,n}\right)$ and $$L_p\left(\Qh;\BC\otimes\cl_{n,n}\right):=L_p\left(\Qh\right)\otimes\left(\BC\otimes \cl_{n,n}\right)~~(1\leq p <\infty)$$
in the following situations (cf.~\cite[Subsection 2.1.b]{hytonen2016analysis} and \cite[Example 2.2.8 of p. 94]{hytonen2016analysis}):
\begin{itemize}
	\item $q=2$ and $1\leq p<\infty$;
	\item $p\leq q\leq 2$.
\end{itemize} 

The systematic treatment of such function spaces properties on the crossroads of pseudo-differential calculus and probability theory will be postponed to forthcoming research papers. At the moment, we can mention the recent paper of Cerejeiras, K\"ahler and Lucas \cite{cerejeiras2021discrete}, which treats several aspects of [weighted] $\ell_q-$ spaces in overlap with symbol classes (see e.g. \cite[Theorem 4 \&  Theorem 8]{cerejeiras2021discrete}).

\subsection{Towards Helmholtz–Leray type decompositions}

In the monograph \cite{gurlebeck1997quaternionic}, G\"urlebeck and Spr\"o\ss ig have popularized the strategy of solving boundary value problems in the context of hypercomplex variables, including the stationary Navier-Stokes equations.
Such strategy, successfully exploited to the discrete setting (cf.~\cite{faustino2006difference,CFV08}), entails the following steps:
\begin{itemize}
	\item[\bf (i)] Compute the fundamental solution of a Dirac type operator;
	\item[\bf (ii)] Determine the right inverse of the Dirac type operator, from the knowledge of the fundamental solution and an analogue of the Borel-Pompeiu formula to incorporate the boundary conditions;
	\item[\bf (iii)] Determine projection operators of Helmholtz–Leray type;
	\item[\bf (iv)] Convert the boundary value problem into an integral equation and solve it, if necessary, by a fixed point scheme.
\end{itemize}

In the context of our framework, one can establish e.g. a parallel with the approach developed in \cite{faustino2006difference}, by considering the pseudo-differential operator $T_h=D_h(-\Delta_h)^{-1}$ -- the so-called {\it discrete Teodorescu operator} (cf.~\cite[p.~239]{gurlebeck1997quaternionic} and \cite[Section 5]{cerejeiras2021discrete}) -- as the right inverse of the discrete Dirac operator $D_h$ ($D_hT_h=I$), as well as $P_h=I-D_h(-\Delta_h)^{-1}D_h$ for the hypercomplex analogue of the so-called {\it Helmholtz–Leray projection}, due to the following set of properties:
\begin{itemize}
	\item[\bf P1] {\it Projection property:} $P_h^2=P_h$ in $\mathcal{S}(h\BZ^n;\BC\otimes \cl_{n,n})$;
	\item[\bf P2] {\it Null property:} $D_hP_h=0$ in $\mathcal{S}(h\BZ^n;\BC\otimes \cl_{n,n})$;
	\item[\bf P3] {\it Direct sum property:} $P_h \Psi=\Psi$, for all $\Psi \in \ker D_h$;
	\item[\bf P4] {\it Multivector functions coming from a potential term :} $P_h(D_h\Psi)=0$, for all $\Psi \in \ker D_h$.
\end{itemize}

Noteworthy, it should be emphasized that such properties seamlessly resembles to the set of properties used in \cite{gurlebeck1997quaternionic} to construct a discrete analogue of the Cauchy integral operator (see \cite[Theorem 5.3.]{gurlebeck1997quaternionic}).

Also, by taking into account the factorization of the semidiscrete Dirac type operators -- considered throughout this paper -- as well as the theory of [fractional semidiscrete] analytic semigroups besides the proof of {\bf Theorem \ref{SemigroupCauchySolver}} and {\bf Theorem \ref{LevyLeblondSolutions}}, the construction of the projection operators depicted on {\bf Figure 3} can be naturally adopted, bearing in mind a possible exploitation of the techniques used e.g. in \cite{CKS05,CFV08}, to obtain hypercomplex formulations of boundary value problems, as well as possible applications on the crossroads of PDEs and stochastics (see e.g. the models considered e.g.~in \cite{KL12,CarvalhoPlanas15,ZP17}). For the technical details, involving the definition of the inverse of the time-fractional/space-fractional operators depicted on the left side of {\bf Figure 3}, we refer to \cite[Section 26.]{SKM93}.

Last but not least, it would be stressed that an important aspect besides the study of L\'evy-Leblond or parabolic type operators -- already quoted on the pioneering work of Cerejeiras, K\"ahler and Sommen (cf.~\cite{CKS05}) -- is the ability to treat non-stationary boundary value problems as stationary ones. Thus, the major reason for including this extra discussion in the end of present paper is mainly for a preliminary motivation for further research studies in the streamlines of this work.

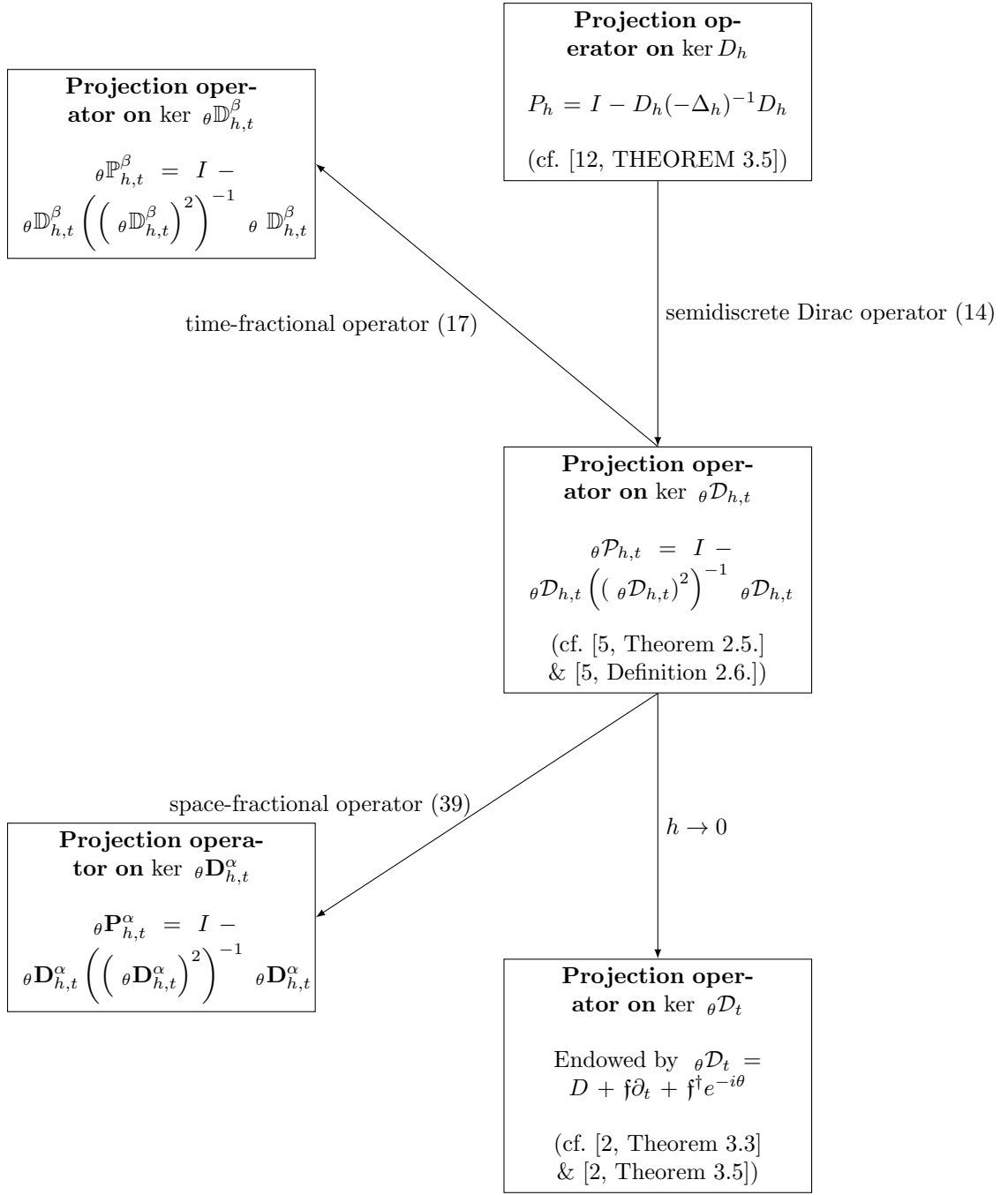
\begin{figure}\label{LerayPicture}
	\begin{tikzpicture}
		\node[mynode] (v1){{\bf Projection operator on $\ker{~}_\theta\mathbb{D}_{h,t}^\beta$} \\ \ \\
			${~}_\theta\mathbb{P}_{h,t}^\beta=I-{~}_\theta\mathbb{D}_{h,t}^\beta\left(\left({~}_\theta\mathbb{D}_{h,t}^\beta\right)^2\right)^{-1}{~}_\theta{~}\mathbb{D}_{h,t}^\beta$};
		\node[mynode,below right=4.0cm of v1](v2) {{\bf Projection operator on $\ker{~}_\theta\mathcal{D}_{h,t}$} \\ \ \\
			${~}_\theta\mathcal{P}_{h,t}=I-{~}_\theta\mathcal{D}_{h,t}\left(\left({~}_\theta\mathcal{D}_{h,t}\right)^2\right)^{-1}{~}_\theta\mathcal{D}_{h,t}$  \\ \ \\
			(cf.~\cite[Theorem 2.5.]{CFV08} \& \cite[Definition 2.6.]{CFV08})};
		
		\node[mynode,below = 8.5cm of v1] (v4){{\bf Projection operator on $\ker{~}_\theta{\bf D}_{h,t}^\alpha$} \\ \ \\
			${~}_\theta{\bf P}_{h,t}^\alpha=I-{~}_\theta{\bf D}_{h,t}^\alpha\left(\left({~}_\theta{\bf D}_{h,t}^\alpha\right)^2\right)^{-1}{~}_\theta{\bf D}_{h,t}^\alpha$};
		\node[mynode,below=4.0cm of v2](v3) {{\bf Projection operator on $\ker{~}_\theta\mathcal{D}_{t}$} \\ \ \\
			Endowed by ${~}_\theta\mathcal{D}_{t}=D+\f\partial_t+\f^\dagger e^{-i\theta}$  \\ \ \\
			(cf.~\cite[Theorem 3.3]{CKS05} \& \cite[Theorem 3.5]{CKS05})};
		\node[mynode, above=4.0cm of v2](v5) {{\bf Projection operator on $\ker D_h$} \\ \ \\
			$P_h=I-D_h(-\Delta_h)^{-1}D_h$  \\ \ \\
			(cf.~\cite[THEOREM 3.5]{faustino2006difference})};
		
		\draw[-latex] (v2.north) -- node[auto,] {{time-fractional operator} (\ref{TimeFractionalDh})} (v1.east);
		\draw[-latex] (v2.south) -- node[left=1.0mm] {{space-fractional operator} (\ref{SpaceFractionalVariant})} (v4.east);
		\draw[-latex] (v5.south) -- node[auto,]{semidiscrete Dirac operator (\ref{SemidiscreteDh})} (v2.north);
		\draw[-latex] (v2.south) -- node[auto,]{$h\rightarrow 0$} (v3.north);
	\end{tikzpicture}
	\caption{The Helmholtz–Leray picture on the fractional semidiscrete case.}
\end{figure}

\bibliography{LevyLeblond_Preprint}

\begin{thebibliography}{10}
\expandafter\ifx\csname url\endcsname\relax
  \def\url#1{\texttt{#1}}\fi
\expandafter\ifx\csname urlprefix\endcsname\relax\def\urlprefix{URL }\fi
\expandafter\ifx\csname href\endcsname\relax
  \def\href#1#2{#2} \def\path#1{#1}\fi

\bibitem{LevyLeblond67}
J.-M. L{\'e}vy-Leblond, Nonrelativistic particles and wave equations,
  Communications in Mathematical Physics 6~(4) (1967) 286--311.

\bibitem{CKS05}
P.~Cerejeiras, U.~K{\"a}hler, F.~Sommen, Parabolic dirac operators and the
  navier--stokes equations over time-varying domains, Mathematical methods in
  the applied sciences 28~(14) (2005) 1715--1724.

\bibitem{CSV07}
P.~Cerejeiras, F.~Sommen, N.~Vieira, Fischer decomposition and special
  solutions for the parabolic dirac operator, Mathematical methods in the
  applied sciences 30~(9) (2007) 1057--1069.

\bibitem{B06}
S.~Bernstein, Factorization of the nonlinear schr{\"o}dinger equation and
  applications, Complex Variables and Elliptic Equations 51~(5-6) (2006)
  429--452.

\bibitem{CFV08}
P.~Cerejeiras, N.~Faustino, N.~Vieira, Numerical clifford analysis for
  nonlinear schr{\"o}dinger problem, Numerical Methods for Partial Differential
  Equations: An International Journal 24~(4) (2008) 1181--1202.

\bibitem{AKTT16}
N.~Aizawa, Z.~Kuznetsova, H.~Tanaka, F.~Toppan, -graded lie symmetries of the
  l{\'e}vy-leblond equations, Progress of Theoretical and Experimental Physics
  2016~(12) (2016) 123A01.

\bibitem{A18}
N.~Aizawa, Generalization of superalgebras to color superalgebras and their
  representations, Advances in Applied Clifford Algebras 28~(1) (2018) 1--14.

\bibitem{FV16}
M.~Ferreira, N.~Vieira, Eigenfunctions and fundamental solutions of the
  fractional laplace and dirac operators: the riemann-liouville case, Complex
  Analysis and Operator Theory 10~(5) (2016) 1081--1100.

\bibitem{FV17}
M.~Ferreira, N.~Vieira, Fundamental solutions of the time fractional
  diffusion-wave and parabolic dirac operators, Journal of Mathematical
  Analysis and Applications 447~(1) (2017) 329--353.

\bibitem{BCBM20}
S.~Bao, D.~Constales, H.~De~Bie, T.~Mertens, Solutions for the l{\'e}vy-leblond
  or parabolic dirac equation and its generalizations, Journal of Mathematical
  Physics 61~(1) (2020) 011509.

\bibitem{BRS20}
D.~Baleanu, J.~E. Restrepo, D.~Suragan, A class of time-fractional dirac type
  operators, Chaos, Solitons \& Fractals 143 (2021) 110590.

\bibitem{faustino2006difference}
N.~Faustino, K.~G{\"u}rlebeck, A.~Hommel, U.~K{\"a}hler, Difference potentials
  for the navier--stokes equations in unbounded domains, Journal of Difference
  Equations and Applications 12~(6) (2006) 577--595.

\bibitem{FK07}
N.~Faustino, U.~Kaehler, Fischer decomposition for difference dirac operators,
  Advances in applied Clifford algebras 17~(1) (2007) 37--58.

\bibitem{FKS07}
N.~Faustino, U.~K{\"a}hler, F.~Sommen, Discrete dirac operators in clifford
  analysis, Advances in Applied Clifford Algebras 17~(3) (2007) 451--467.

\bibitem{BSSV08}
F.~Brackx, H.~De~Schepper, F.~Sommen, L.~Van~de Voorde, Discrete clifford
  analysis: a germ of function theory, in: Hypercomplex Analysis, Springer,
  2008, pp. 37--53.

\bibitem{RSKS10}
H.~De~Ridder, H.~De~Schepper, U.~K{\"a}hler, F.~Sommen, Discrete function
  theory based on skew weyl relations, Proceedings of the American Mathematical
  Society 138~(9) (2010) 3241--3256.

\bibitem{KS75}
J.~Kogut, L.~Susskind, Hamiltonian formulation of wilson's lattice gauge
  theories, Physical Review D 11~(2) (1975) 395.

\bibitem{Rabin82}
J.~M. Rabin, Homology theory of lattice fermion doubling, Nuclear Physics B
  201~(2) (1982) 315--332.

\bibitem{KK04}
I.~Kanamori, N.~Kawamoto, Dirac--kaehler fermion from clifford product with
  noncommutative differential form on a lattice, International Journal of
  Modern Physics A 19~(05) (2004) 695--736.

\bibitem{Sushch14}
V.~Sushch, A discrete model of the dirac-k{\"a}hler equation, Reports on
  Mathematical Physics 73~(1) (2014) 109--125.

\bibitem{gurlebeck1997quaternionic}
K.~G{\"u}rlebeck, W.~Spr{\"o}ssig, Quaternionic and Clifford calculus for
  physicists and engineers, Vol.~1, John Wiley \& Sons, 1997.

\bibitem{VR16}
J.~Vaz~Jr, R.~da~Rocha~Jr, An introduction to Clifford algebras and spinors,
  Oxford University Press, 2016.

\bibitem{FaustinoKGordonDirac16}
N.~Faustino, Solutions for the klein--gordon and dirac equations on the lattice
  based on chebyshev polynomials, Complex analysis and operator theory 10~(2)
  (2016) 379--399.

\bibitem{FaustinoMMAS17}
N.~J. Rodrigues~Faustino, A conformal group approach to the dirac--k{\"a}hler
  system on the lattice, Mathematical Methods in the Applied Sciences 40~(11)
  (2017) 4118--4127.

\bibitem{F19b}
N.~Faustino, A note on the discrete cauchy-kovalevskaya extension, Mathematical
  Methods in the Applied Sciences 42~(4) (2019) 1312--1320.

\bibitem{LizamaRoncal18}
C.~Lizama, L.~Roncal, H{\"o}lder-lebesgue regularity and almost periodicity for
  semidiscrete equations with a fractional laplacian (2018).

\bibitem{F19}
N.~Faustino, Relativistic wave equations on the lattice: an operational
  perspective, in: Topics in Clifford Analysis, Springer, 2019, pp. 439--469.

\bibitem{SKM93}
S.~G. Samko, A.~A. Kilbas, O.~I. Marichev, Fractional integrals and
  derivatives: theory and applications (1993).

\bibitem{PBM86}
A.~Prudnikov, Y.~A. Brychkov, O.~I. Marichev, Integrals and series, volume 1:
  Elementary functions, Gordon\&Breach Sci. Publ., New York (1986).

\bibitem{CiaurriGRTV18}
O.~Ciaurri, L.~Roncal, P.~R. Stinga, J.~L. Torrea, J.~L. Varona, Nonlocal
  discrete diffusion equations and the fractional discrete laplacian,
  regularity and applications, Advances in Mathematics 330 (2018) 688--738.

\bibitem{GLizamaM21}
J.~Gonz{\'a}lez-Camus, C.~Lizama, P.~J. Miana, Fundamental solutions for
  semidiscrete evolution equations via banach algebras, Advances in Difference
  Equations 2021~(1) (2021) 1--32.

\bibitem{CLRV15}
{\'O}.~Ciaurri, C.~Lizama, L.~Roncal, J.~L. Varona, On a connection between the
  discrete fractional laplacian and superdiffusion, Applied Mathematics Letters
  49 (2015) 119--125.

\bibitem{KL02}
V.~V. Kulish, J.~L. Lage, Application of fractional calculus to fluid
  mechanics, J. Fluids Eng. 124~(3) (2002) 803--806.

\bibitem{KL12}
V.~Keyantuo, C.~Lizama, On a connection between powers of operators and
  fractional cauchy problems, Journal of Evolution Equations 12~(2) (2012)
  245--265.

\bibitem{BMN09}
B.~Baeumer, M.~Meerschaert, E.~Nane, Brownian subordinators and fractional
  cauchy problems, Transactions of the American Mathematical Society 361~(7)
  (2009) 3915--3930.

\bibitem{MS12}
M.~M. Meerschaert, A.~Sikorskii, Stochastic models for fractional calculus, de
  Gruyter, 2019.

\bibitem{BBRS14}
F.~Baaske, S.~Bernstein, H.~De~Ridder, F.~Sommen, On solutions of a discretized
  heat equation in discrete clifford analysis, Journal of Difference Equations
  and Applications 20~(2) (2014) 271--295.

\bibitem{Mainardi01}
F.~Mainardi, Y.~Luchko, G.~Pagnini, The fundamental solution of, Fractional
  Calculus and Applied Analysis 4~(2) (2001) 153--192.

\bibitem{F20}
N.~Faustino, Time-changed dirac--fokker--planck equations on the lattice,
  Journal of Fourier Analysis and Applications 26 (2020) 1--31.

\bibitem{B16}
S.~Bernstein, A fractional dirac operator, in: Noncommutative Analysis,
  Operator Theory and Applications, Springer, 2016, pp. 27--41.

\bibitem{CarvalhoPlanas15}
P.~M. de~Carvalho-Neto, G.~Planas, Mild solutions to the time fractional
  navier--stokes equations in rn, Journal of Differential Equations 259~(7)
  (2015) 2948--2980.

\bibitem{ZP17}
Y.~Zhou, L.~Peng, On the time-fractional navier--stokes equations, Computers \&
  Mathematics with Applications 73~(6) (2017) 874--891.

\bibitem{RuzhanskyT10}
M.~Ruzhansky, V.~Turunen, Pseudo-differential operators and symmetries:
  background analysis and advanced topics, Vol.~2, Springer Science \& Business
  Media, 2009.

\bibitem{CiaurriGRTV17}
{\'O}.~Ciaurri, T.~A. Gillespie, L.~Roncal, J.~L. Torrea, J.~L. Varona,
  Harmonic analysis associated with a discrete laplacian, Journal d'Analyse
  Math{\'e}matique 132~(1) (2017) 109--131.

\bibitem{hytonen2016analysis}
T.~Hyt{\"o}nen, J.~Van~Neerven, M.~Veraar, L.~Weis, Analysis in Banach spaces:
  Volume I: Martingales and Littlewood-Paley Theory, Vol.~63, Springer, 2016.

\bibitem{cerejeiras2021discrete}
P.~Cerejeiras, U.~Kaehler, S.~Lucas, Discrete pseudo-differential operators in
  hypercomplex analysis, Mathematical Methods in the Applied Sciences (2021).

\end{thebibliography}

\end{document}